\documentclass[oneside, reqno]{amsart}
\pdfoutput=1
\usepackage{subfiles}
\usepackage[utf8]{inputenc}
\usepackage[T1]{fontenc}
\usepackage[dvipsnames]{xcolor}
\usepackage{ifdraft}
\usepackage{amssymb}
\usepackage{amsmath}
\usepackage{amsthm}
\usepackage{caption} 
\usepackage{microtype}
\usepackage{mathtools}
\usepackage{tikz}
\usepackage[bookmarks, bookmarksdepth=subsection, colorlinks, allcolors=Blue, final]{hyperref}
\usepackage{tikz-cd}
\usepackage{lipsum}
\usepackage[foot]{amsaddr}

\ifdraft{
    \usepackage{showlabels}
}{}

\usepackage[colorinlistoftodos, obeyDraft, textsize=tiny]{todonotes}

\usepackage[shortlabels]{enumitem}
\setlist[itemize]{leftmargin=2\parindent}
\setlist{topsep=0ex, itemsep=0.4ex, parsep=0ex}
\setenumerate[0]{label=(\roman*), font=\normalfont}

\usepackage[capitalize, noabbrev, nameinlink]{cleveref}
\RemoveFromHook{label}[firstaid/cleveref]
\Crefname{equation}{}{}
\Crefname{enumi}{}{}
\Crefname{axiom}{Axiom}{Axioms}
\Crefname{diagram}{Diagram}{Diagrams}
\crefformat{footnote}{#2\footnotemark[#1]#3}

\usepackage[
    backend=biber,
    style=numeric-comp,
    giveninits=true,
    maxbibnames=10,
    minbibnames=10,
    mincrossrefs=10,
    bibencoding=utf8]{biblatex} 

\AtBeginBibliography{
    \raggedright
}

\ifSubfilesClassLoaded{
    \AtEndDocument{%
        \printbibliography
    }
}{}

\makeatletter
\renewcommand{\@mkboth}[2]{}
\makeatother

\usetikzlibrary{decorations.markings}
\usetikzlibrary{cd}

\tikzset{
    ulcorner/.style={
        label={[inner sep=0pt, anchor=center]south east:{\lrcorner}}
    },
    urcorner/.style={
        label={[inner sep=0pt, anchor=center]south west:{\llcorner}}
    },
    llcorner/.style={
        label={[inner sep=0pt, anchor=center]north east:{\urcorner}}
    },
    lrcorner/.style={
        label={[inner sep=0pt, anchor=center]north west:{\ulcorner}}
    },
}

\tikzcdset{
    diagrams={ampersand replacement=\&},
    postmark/.style={
        postaction={
            decorate,
            decoration={
                markings,
                mark={#1}
            }
        }
    },
    display style/.style={
        cells={font=\everymath\expandafter{\the\everymath\displaystyle}},
    },
    dagger arrow/.style={
        postmark=at position 0.5 with {
            \draw[-] (#1,1.5pt) -- (#1,-1.5pt);
        },
        labels={
            inner sep=0.75ex,
        }
    },
    mono/.style={tail},
    epi/.style={two heads},
    dagger mono/.style={dagger arrow={0pt}, mono},
    dagger epi/.style={dagger arrow={-1pt}, epi},
    left adjoint/.style={
        phantom,
        "\dashv"{font=\scriptsize},
        sloped,
        rotate={90}
    },
    iso/.style={"\cong"{sloped}},
}
\input{preamble/widebar}
\theoremstyle{plain}
\newtheorem{theorem}                {Theorem}     [section]

\newtheorem{lemma}        [theorem] {Lemma}
\newtheorem{proposition}  [theorem] {Proposition}
\newtheorem{corollary}    [theorem] {Corollary}
\newtheorem*{theorem*}              {Theorem}

\theoremstyle{definition}

\newtheorem{definition}   [theorem] {Definition}

\newtheorem*{definition*} {Definition}

\theoremstyle{remark}
\newtheorem{remark}       [theorem] {Remark}

\newtheorem*{remark*}               {Remark}

\renewcommand{\leq}{\leqslant}

\renewcommand{\geq}{\geqslant}

\newcommand{\slt}{\prec}
\newcommand{\sgt}{\succ}

\newcommand{\hquad}{\mkern9mu}

\newcommand{\blank}{\,\underline{\phantom{n}}\,}

\DeclarePairedDelimiter{\paren}{(}{)}
\DeclarePairedDelimiter{\abs}{|}{|}
\DeclarePairedDelimiter{\norm}{\|}{\|}
\DeclarePairedDelimiterXPP{\innerProd}[2]{\bgroup}{\langle}{\rangle}{\egroup}{#1 \delimsize\vert\mathopen{} #2}
\DeclarePairedDelimiterX{\setb}[2]{\lbrace}{\rbrace}{#1 \,\delimsize\vert\,\mathopen{} #2}
\DeclarePairedDelimiter{\set}{\lbrace}{\rbrace}

\newcommand{\zero}{\mathrm{O}}

\newcommand{\Pos}[1]{{#1}_{\geq 0}}
\newcommand{\SelfAdj}[1]{{#1}_{\mathrm{H}}}
\newcommand{\inv}{\ast}

\newcommand{\field}[1]{\mathbb{#1}}
\newcommand{\Nats}{\field{N}}
\newcommand{\Ints}{\field{Z}}
\newcommand{\Rats}{\field{Q}}

\newcommand{\Reals}{\field{R}}
\newcommand{\PosReals}{\Pos{\Reals}}
\newcommand{\Comps}{\field{C}}
\newcommand{\Quats}{\field{H}}

\newcommand{\cat}[1]{\mathbf{#1}}
\newcommand{\pcat}[2]{\cat{#1}_{#2}}

\newcommand{\Con}[1]{{#1}_{\leq 1}}
\newcommand{\NMono}[1]{{#1}_{1}}
\newcommand{\Hilb}{\cat{Hilb}}

\newcommand{\PHilb}[1]{\pcat{Hilb}{#1}}
\newcommand{\Mat}[1]{\pcat{Mat}{#1}}
\newcommand{\WeightMat}[1]{\pcat{WeightMat}{#1}}
\newcommand{\InnerProd}[1]{\pcat{InnProd}{#1}}
\newcommand{\Mod}[1]{\pcat{Mod}{#1}}
\newcommand{\FinInnerProd}[1]{\pcat{FInnProd}{#1}}

\newcommand{\Ab}{\cat{Ab}}
\newcommand{\Set}{\cat{Set}}

\newcommand{\C}{\cat{C}}

\newcommand{\D}{\cat{D}}

\newcommand{\pair}[2]{\begin{bsmallmatrix} #1 \\ #2 \end{bsmallmatrix}}
\newcommand{\pairBig}[2]{\begin{bmatrix} #1 \\ #2 \end{bmatrix}}
\newcommand{\copair}[2]{\begin{bsmallmatrix} #1 & #2 \end{bsmallmatrix}}
\newcommand{\copairBig}[2]{\begin{bmatrix} #1 & #2 \end{bmatrix}}

\DeclareMathOperator{\Ker}{Ker}
\DeclareMathOperator{\Eq}{Eq}

\newcommand{\Ran}{\operatorname{Ran}}

\renewcommand{\Re}{\operatorname{Re}}
\renewcommand{\Im}{\operatorname{Im}}

\DeclareMathOperator{\centre}{Z}
\newcommand{\diagonal}{\Delta}
\newcommand{\codiagonal}{\nabla}

\newcommand{\bsmallmat}[1]{\begin{bsmallmatrix}#1\end{bsmallmatrix}}
\makeatletter
\newcommand{\fps@diagram}{tbp}
\newcounter{diagram}
\def\ftype@diagram{1}
\def\ext@diagram{lof}
\def\fnum@diagram{\diagramname\ \thediagram}
\def\diagramname{Diagram}
\newenvironment{diagram}{%
\@float{diagram}%
}{%
\end@float
}
\newenvironment{diagram*}{%
\@dblfloat{diagram}%
}{%
\end@dblfloat
}
\makeatother

\usepackage{bm}

\reversemarginpar

\begin{document}
\title[Pre-Hilbert \(*\)-categories]{Pre-Hilbert \texorpdfstring{\(\bm{*}\)}{*}-categories\\[2ex]\footnotesize The Hilbert-space analogue of\\abelian categories}

\date{\today}
\author{Matthew {Di Meglio}}
\address{The University of Edinburgh, 10 Crichton St, Edinburgh, United Kingdom}
\email{m.dimeglio@ed.ac.uk}
\subjclass[2020]{18M40, 18E99, 46M15, 06F25}
\keywords{Hilbert space, category, order, dilator}

\begin{abstract}
    This article introduces pre\nobreakdash-Hilbert \(\inv\)\nobreakdash-categories---an abstraction of categories exhibiting ``algebraic'' aspects of the theory of Hilbert spaces. Notably, finite biproducts in pre\nobreakdash-Hilbert \(\inv\)\nobreakdash-categories can be orthogonalised using the Gram--Schmidt process, and generalised notions of positivity and contraction support variants of Douglas' lemma and Sz.-Nagy's unitary dilation theorem. Underpinning these generalisations is the structure of an involutive identity-on-objects contravariant endofunctor, which encodes adjoints of morphisms. The axioms for pre\nobreakdash-Hilbert \(\inv\)\nobreakdash-categories are otherwise inspired by those for abelian categories, comprising a few simple properties of products and kernels. Additivity is not assumed, but nevertheless follows. In fact, the similarity with abelian categories runs deeper---pre\nobreakdash-Hilbert \(\inv\)\nobreakdash-categories are quasi-abelian and thus also homological. Examples include the \(\inv\)\nobreakdash-category of unitary representations of a groupoid, the \(\inv\)\nobreakdash-category of finite-dimensional inner product modules over an ordered division \(\inv\)\nobreakdash-ring, and the \(\inv\)\nobreakdash-category of self-dual Hilbert modules over a W*-algebra.
\end{abstract}
\maketitle

\setcounter{tocdepth}{1}
\tableofcontents

\section{Introduction}
Surprisingly, the categories of real, complex and quaternionic Hilbert spaces and bounded linear maps, respectively denoted \(\Hilb_\Reals\), \(\Hilb_\Comps\) and \(\Hilb_\Quats\), can be characterised among all categories by a few simple categorical axioms~\cite{heunen:2022:axioms-for-the-category,tobin:characterisations-category-hilbert,tobin-lack:hilb}. Remarkably, this characterisation does not refer to any analytic or topological notions like real numbers, limits of sequences, norms, continuity or probability.

The foundation of this characterisation is the notion of \(\inv\)\nobreakdash-category. A \textit{\(\inv\)\nobreakdash-category} is a category equipped with a functorial and involutive choice of morphisms \(f^\inv \colon Y \to X\) for each morphism \(f \colon X \to Y\).\footnote{The name \textit{dagger category} for this concept and the symbol \(\dagger\) for the involution was popularised by the school of categorical quantum mechanics~\cite{abramsky:2004:categorical-semantics-quantum-protocols,selinger:2007:dagger-compact-closed-categories,harding:2009:link-between-quantum-logic}. This article adopts the name \textit{\(\inv\)\nobreakdash-category} for consistency with the literature on \(\inv\)\nobreakdash-rings and C*-algebras. In this article, \(\inv\)\nobreakdash-categories are not, \textit{a priori}, enriched in the category of complex vector spaces. We reserve the term \textit{\(\Comps\)-linear \(\inv\)\nobreakdash-category} for such an enriched \(\inv\)\nobreakdash-category.} The operation \(f \mapsto f^\inv\) will be referred to as the \textit{involution} of the \(\inv\)\nobreakdash-category. An involution on a category encodes the same kind of geometric information as an inner product on a vector space.

The rest of the characterisation comprises several axioms about products and equalisers and their compatibility with the involution. Concretely, it says that a \(\inv\)\nobreakdash-category is equivalent to \(\Hilb_{\field{K}}\) for some \(\field{K}\) among \(\Reals\), \(\Comps\) and \(\Quats\) if and only if
\begin{enumerate}[(H1)]
    \item
    there is a zero object,
    \item
    every pair of objects has an orthonormal biproduct,\footnote{Isometry, isometric equaliser, and orthonormal biproduct are respectively called \textit{dagger monomorphism}, \textit{dagger equaliser}, and \textit{dagger biproduct} in the literature on dagger categories.\label{dagger-names}}
    \item
    every parallel pair of morphisms has an isometric equaliser,\cref{dagger-names}
    \item
    every isometry is a normal monomorphism (i.e., a kernel),\cref{dagger-names}
    \item 
    \label{axiom:hilbert-directed}
    the wide subcategory of isometries has directed colimits, and
    \item there is a simple separating object.
\end{enumerate}
Elementary concepts from \(\inv\)\nobreakdash-category theory, including \textit{isometry} and \textit{orthonormal biproduct}, are recalled in \cref{sect:background}. The reader not yet familiar with them should, at this point, still appreciate the simplicity and elegance of the characterisation.

There is a strikingly similar well-known characterisation of the categories \(\Mod{R}\) of right modules over a ring \(R\). One way to phrase it\footnote{In modern terminology, \textcite[106]{freyd:1964:abelian-categories} states (without proof) that a category is equivalent to the category \(\Mod{R}\) for some ring \(R\) if and only if it is abelian, cocomplete, and has a small-projective separating object. Axioms \cref*{axiom:abelian-zero,axiom:abelian-product,axiom:abelian-kernel,axiom:abelian-normal} form \citeauthor{borceux:1994:handbook-categorical-algebra-2}'s definition~\cite[Definition~1.4.1]{borceux:1994:handbook-categorical-algebra-2} of abelian category. Recall that cocompleteness is equivalent to the existence of coequalisers and arbitrary coproducts, and that abelian categories already have coequalisers and finite biproducts. With axiom \cref*{axiom:abelian-directed}, arbitrary coproducts can be constructed from finite ones.} is that a category is equivalent to the category \(\Mod{R}\) for some ring~\(R\) if and only if
\begin{enumerate}[label={(A\arabic*)}]
    \item 
    \label{axiom:abelian-zero}
    there is a zero object,
    \item
    \label{axiom:abelian-product}
    every pair of objects has both a product and a coproduct,
    \item
    \label{axiom:abelian-kernel}
    every morphism has both a kernel and a cokernel,
    \item
    \label{axiom:abelian-normal}
    every monomorphism is normal and every epimorphism is normal,
    \item 
    \label{axiom:abelian-directed}
    every directed diagram of split monomorphisms has a colimit, and
    \item
    \label{axiom:abelian-separator}
    there is a small-projective separating object.
\end{enumerate}
Written in this way, the axioms for \(\Mod{R}\) correspond one-to-one with the axioms for \(\Hilb_\field{K}\). In fact, this correspondence is stronger than it might first appear: the duals of the axioms for \(\Hilb_\field{K}\) hold automatically because \(\inv\)\nobreakdash-categories are self-dual.

Axioms \cref*{axiom:abelian-zero,axiom:abelian-product,axiom:abelian-kernel,axiom:abelian-normal} form one~\cite[Definition~1.4.1]{borceux:1994:handbook-categorical-algebra-2} of the equivalent definitions of \textit{abelian category}. The notion of abelian category is a well-studied abstraction of the category of abelian groups, with important applications across many branches of pure mathematics, including algebra, algebraic geometry, and topology. This raises the central question in this article: what can be said about the \(\inv\)\nobreakdash-category analogue of abelian category, that is, the \(\inv\)\nobreakdash-categories that satisfy axioms \cref*{axiom:hilbert-zero,axiom:hilbert-product,axiom:hilbert-equaliser,axiom:hilbert-normal}?

\begin{definition*}
A \textit{pre\nobreakdash-Hilbert \(\inv\)\nobreakdash-category}\footnote{
Unfortunately, the name \textit{pre\nobreakdash-Hilbert category} is already used by \textcite{heunen:2009:hilbert-categories} for the symmetric monoidal variant of the concept defined here. \citeauthor{heunen:2009:hilbert-categories} also strived for a \(\inv\)\nobreakdash-category analogue of abelian categories, but this analogy was marred by the reliance on the monoidal product to prove abelian-like properties, such as additivity, that ought to hold without it. Pre\nobreakdash-Hilbert \(\inv\)\nobreakdash-categories, as defined in this article, are closer in spirit to abelian categories, and also include examples that do not have a tensor product, such as \(\Hilb_\Quats\). Given that the tensor product of Hilbert spaces is central to categorical quantum mechanics~\cite{abramsky:2004:categorical-semantics-quantum-protocols,harding:2009:link-between-quantum-logic,selinger:2007:dagger-compact-closed-categories}, the fact that it is not actually essential to the category theory of Hilbert spaces is quite a revelation.} is a \(\inv\)\nobreakdash-category in which
\begin{enumerate}[(R1)]
    \item
    there is a zero object,
    \item
    every pair of objects has an orthonormal biproduct,
    \item
    every morphism has an isometric kernel,
    \item
    every diagonal \(\diagonal \colon X \to X \oplus X\) is a normal monomorphism.
\end{enumerate}
\end{definition*}

The careful reader will have noticed that axioms \cref*{axiom:R-zero,axiom:R-coproduct,axiom:R-kernel,axiom:R-normal} are slightly different from axioms \cref*{axiom:hilbert-zero,axiom:hilbert-product,axiom:hilbert-equaliser,axiom:hilbert-normal}. These two families of axioms are actually equivalent to each other (\cref{prop:equivalent-definition-heunen}), but axioms \cref*{axiom:R-zero,axiom:R-coproduct,axiom:R-kernel,axiom:R-normal} align even more closely with the corresponding axioms for abelian categories.

The parallels between pre\nobreakdash-Hilbert \(\inv\)\nobreakdash-categories and abelian categories run deeper than merely the axioms. Firstly, pre\nobreakdash-Hilbert \(\inv\)\nobreakdash-categories are additive. Ensuring that the canonical addition of morphisms (which comes from the biproducts) admits additive inverses is, in fact, the purpose of axiom \cref{axiom:R-normal}. That this can be done with an axiom that so faithfully parallels the corresponding axiom \cref{axiom:abelian-normal} for abelian categories, and yet also holds in \(\Hilb_{\Comps}\) (unlike axiom \cref{axiom:abelian-normal}), was not at all obvious. With axiom \cref{axiom:R-normal}, additivity is a corollary of the following novel\footnote{
\cref{prop:semiadditive-abelian-objects} is, to the best of the author's knowledge, new. While it is reminiscent of several known characterisations of abelian objects (see, e.g.,~\cite[Proposition~9]{bourn:normal-subobjects-abelian} and~\cite[Proposition~1.5.11]{borceux-bourn:2004:malcev-protomodular-homological}), those results do not apply to all categories with finite biproducts.
} characterisation (\cref{prop:semiadditive-abelian-objects}) of \textit{abelian objects}---objects whose identity morphism has an additive inverse---in categories with finite biproducts.

\begin{theorem*}
An object \(X\) in a category with finite biproducts is abelian if and only if the diagonal \(\diagonal \colon X \to X\oplus X\) is a kernel of a split epimorphism. 
\end{theorem*}

\noindent Secondly, in pre\nobreakdash-Hilbert \(\inv\)\nobreakdash-categories, normal monomorphisms are stable under pushout along arbitrary morphisms (\cref{prop:kernels-pushout-stable}). Hence pre\nobreakdash-Hilbert \(\inv\)\nobreakdash-categories are quasi-abelian (see \cite{schneiders:1999:quasi-abelian-sheaves,rump:2001:almost-abelian}), and thus also homological (see \cite[Remarks~5.3]{rosicky-tholen:2007:factorisation-fibration}). From this follows many useful properties of pre\nobreakdash-Hilbert \(\inv\)\nobreakdash-categories, such as the diagram lemmas from homological algebra. These similarities between pre\nobreakdash-Hilbert \(\inv\)\nobreakdash-categories and abelian categories are the theme of \cref{sect:abelian-properties}.

While the notion of orthogonal complement makes sense in any \(\inv\)\nobreakdash-category with isometric kernels~\cite{heunen:2010:quantum-logic-dagger-kernel,schwab:2012:quantum-logic-dagger-kernel}, it is especially well behaved in pre\nobreakdash-Hilbert \(\inv\)\nobreakdash-categories because of additivity. For example, finite biproducts in pre\nobreakdash-Hilbert \(\inv\)\nobreakdash-categories can be transformed into \textit{orthogonal biproducts} (\cref{def:orthogonal-biproduct}) using a generalisation of the Gram--Schmidt process (\cref{prop:gram-schmidt}). These aspects of orthogonality are the topic of \cref{sect:gram-schmidt}.

Just as the category \(\Ab = \Mod{\Ints}\) of abelian groups is the prototypical abelian category, the \(\inv\)\nobreakdash-categories \(\Hilb_\Reals\), \(\Hilb_\Comps\) and \(\Hilb_\Quats\) of Hilbert spaces (which are not abelian) are the prototypical pre\nobreakdash-Hilbert \(\inv\)\nobreakdash-categories. Other examples include:
\begin{itemize}
    \item for each groupoid \(G\), the \(\inv\)\nobreakdash-category of unitary representations of \(G\) and their intertwiners (\cref{subsect:r-cat});
    \item for each ordered division \(\inv\)\nobreakdash-ring \(D\), the \(\inv\)\nobreakdash-category of finite-dimensional inner product \(D\)\nobreakdash-modules and \(D\)-linear maps (\cref{sect:fin-inner-prod-modules}); and
    \item for each complex W*-algebra \(A\), the \(\inv\)\nobreakdash-category of self-dual Hilbert \(A\)\nobreakdash-modules and bounded \(A\)-linear maps (\cref{sect:hilbert-W-modules}).
\end{itemize}
The common theme of \textit{inner products} among these examples is no coincidence. Given a pre\nobreakdash-Hilbert \(\inv\)\nobreakdash-category \(\C\), each hom-set \(\C(A, A)\) is canonically a (partially) ordered \(\inv\)\nobreakdash-ring, and each hom-set \(\C(A, X)\) is canonically an inner product module over \(\C(A, A)\), as is explained in \cref{prop:canonical-structure}. In brief, the canonical partial ordering of the Hermitian elements of \(\C(A,A)\) is defined by \(a \leq b\) if \(b - a = y^\inv y\) for some object \(Y\) and some element \(y \in \C(A, Y)\), and the canonical \(\C(A,A)\)-valued inner product on \(\C(A, X)\) is defined by \(\innerProd{x_1}{x_2} = {x_1}^{\!\inv} x_2\). These connections with \(\inv\)\nobreakdash-categories of inner product modules are the topic of \cref{sect:examples}.

The partial order defined above exhibits the following unexpected and useful property: if \(a \geq 1\), then \(a\) is invertible (\cref{prop:symmetric}). In unbounded operator algebra, this property is called \textit{symmetry} (see, e.g., \cite{schotz:equivalence-order-algebraic,schmudgen:2012:unbounded-operators,dixon:generalized-b-algebras,inoue:class-unbounded-operator}). It is a useful tool for extending the theory of bounded operators to unbounded operators. Using symmetry, arbitrary morphisms of a pre\nobreakdash-Hilbert \(\inv\)\nobreakdash-category can be written as the composite of a contraction and the inverse of a contraction (see below). Also, pre\nobreakdash-Hilbert \(\inv\)\nobreakdash-categories are uniquely \(\Rats\)\nobreakdash-linear (\cref{prop:rat-linear}), that is, enriched in the category of rational vector spaces. These connections between positivity and inverses are discussed in \cref{sect:positivity-inverses}.

Douglas' lemma~\cite[Theorem~1]{douglas:majorization-factorization-range} is a well-known and useful result identifying range inclusion, majorisation, and factorisation of operators on a Hilbert space. Its main claim is that the following statements about bounded linear maps \(f \colon X \to Y\) and \(y \colon A \to Y\) between Hilbert spaces are equivalent:
\begin{enumerate}
    \item \(\Ran y \subseteq \Ran f\),
    \item \(yy^\inv \leq rff^\inv\) for some \(r \in \Reals_{\geq 0}\), and
    \item \(y = fx\) for some bounded linear operator \(x \colon A \to X\).
\end{enumerate}
Remarkably, almost all of Douglas' lemma holds more generally in every pre\nobreakdash-Hilbert \(\inv\)\nobreakdash-category. Proving this is the focus of \cref{s:douglas}.

The familiar notion of contraction from the theory of Hilbert spaces has the following generalisation: a morphism \(f\) in a pre\nobreakdash-Hilbert \(\inv\)\nobreakdash-category is a \textit{contraction} if \(f^\inv \! f \leq 1\). These contractions form a wide \(\inv\)\nobreakdash-subcategory (\cref{prop:con-is-subcategory}) with properties similar to those of the \(\inv\)\nobreakdash-categories of Hilbert spaces and contractions (\cref{prop:contraction-properties}). This is the subject of \cref{sect:category-con}.

Sz.-Nagy's unitary dilation theorem is the foundation of the modern theory of contractions on Hilbert spaces~\cite{nagy:2010:harmonic-analysis}. It says that for each contraction \(f \colon X \to X\) on a Hilbert space \(X\), there is an isometry \(s \colon X \to S\) and a unitary operator \(u \colon S \to S\) such that \(u\) is a \textit{dilation} of \(f\), that is, \(f^n = s^\inv u^n s\) for all positive integers \(n\), and such that \(S\) is \textit{minimal} in a certain sense. For category theorists, it is natural to seek a universal property capturing minimality, as well as a notion of minimal dilation for contractions between different Hilbert spaces. These goals lead to the following new universal construction, which is introduced in \cref{s:codilator}.

\begin{definition*}
A \textit{codilation} of a morphism \(f \colon X \to Y\) in a \(\inv\)\nobreakdash-category is a cospan \((T, t_1, t_2)\) where \(t_1 \colon X \to T\) and \(t_2 \colon Y \to T\) are isometries and \({t_2}^{\!\inv} t_1 = f\). A \textit{codilator} of \(f\) is a codilation \((S, s_1, s_2)\) of \(f\) such that, for all codilations \((T, t_1, t_2)\) of \(f\), there is a unique isometry \(t \colon S \to T\) such that \(ts_1 = t_1\) and \(ts_2 = t_2\), as depicted in \cref{f:codilator}.
\begin{diagram}
    \centering
    \begin{tikzcd}[sep=large]
        \&
    T
        \&
    \\
        \&
    S
        \arrow[u, "t"{pos=0.3}]
        \&
    \\[-2em]
    X
        \arrow[uur, "t_1"]
        \arrow[ur, "s_1" swap]
        \&\&
    Y
        \arrow[uul, "t_2" swap]
        \arrow[ul, "s_2"]
    \end{tikzcd}
    \caption{}
    \label{f:codilator}
\end{diagram}
\end{definition*}

The following variant of Sz.-Nagy's theorem is proved in \cref{s:codilator}.

\begin{theorem*}
Every contraction in a pre\nobreakdash-Hilbert \(\inv\)\nobreakdash-category has a codilator.
\end{theorem*}

The notion of codilator transcends pre\nobreakdash-Hilbert \(\inv\)\nobreakdash-categories. For example, the category of fully supported finite probability spaces and stochastic maps is a \(\inv\)\nobreakdash-category where the involution is Bayesian inversion, and each morphism in this \(\inv\)\nobreakdash-category has a codilator, namely, its \textit{bloom-shriek factorisation}~\cite{fullwood:information-stochastic}.\footnote{To be pedantic, the codilator of such a stochastic map is the cospan formed by its \textit{bloom} and the Bayesian inverse of its \textit{shriek}.} Additionally, independent pullbacks~\cite{simpson:category-theoretic-structure-independence,simpson:equivalence-conditional-independence} are related to codilators of partial isometries. An account of codilators in arbitrary \(\inv\)\nobreakdash-categories will be given in a future article.

Pre\nobreakdash-Hilbert \(\inv\)\nobreakdash-categories that also satisfy axiom~\cref*{axiom:hilbert-directed}, that is, those whose wide subcategory of isometries has directed colimits, will be called \textit{Hilbert \(\inv\)\nobreakdash-categories}. Examples include the \(\inv\)\nobreakdash-category of unitary representations of a groupoid and the \(\inv\)\nobreakdash-category of self-dual Hilbert modules over a complex W*-algebra. For each Hilbert \(\inv\)\nobreakdash-category~\(\C\) and all objects \(A\) and \(X\) of \(\C\), the ordered \(\inv\)\nobreakdash-rings \(\C(A, A)\) have suprema of bounded increasing nets, and the inner product modules \(\C(A, X)\) are also complete in a certain sense.
Building on the ideas of \textcite{dimeglio:2024:con}, both kinds of (analytic) limits of nets of elements may be extracted from (categorical) colimits of associated directed diagrams of contractions. The fact that every contraction has a codilator is the key to proving that the wide subcategory of contractions of a Hilbert \(\inv\)\nobreakdash-category has directed colimits. The theory of Hilbert \(\inv\)\nobreakdash-categories will be expounded in a sequel to this article.

\subsection*{Acknowledgements}
Many thanks go to Chris Heunen, my PhD supervisor, for his help and support throughout the development of this work. Thanks also go to Nadja Egner, Marino Gran and Tim Van der Linden for pointing out the connection between pre\nobreakdash-Hilbert \(\inv\)\nobreakdash-categories and quasi-abelian categories, and to Paolo Perrone for pointing out the connection between codilators and the bloom-shriek factorisation. I am also very grateful to Steve Lack, Tom Leinster, and the anonymous reviewers for their helpful suggestions and corrections. \Cref{p:fact-eq} was proved in collaboration with GPT-5.6 Sol.

\section{Understanding the definition}
\label{sect:background}

This section introduces the concepts, terminology and notation from \(\inv\)\nobreakdash-category theory needed to understand the definition of pre\nobreakdash-Hilbert \(\inv\)\nobreakdash-category, which is given at the end, in \cref{subsect:r-cat}. The concept of \textit{closed monomorphism} (\cref{def:closed-mono}) is new. Everything else before \cref{subsect:r-cat} is well known~\cite{borceux:1994:handbook-categorical-algebra-2,heunen:2019:categories-for-quantum-theory,heunen:2010:quantum-logic-dagger-kernel,vicary:2011:complex-numbers}, and so is presented without proof. As many concepts in \(\inv\)\nobreakdash-category theory are lifted from the theory of Hilbert spaces, the \(\inv\)\nobreakdash-category \(\Hilb\) of complex Hilbert spaces and bounded linear maps is a running example. More examples are described in \cref{subsect:r-cat} and \cref{sect:examples}.

\subsection{Hilbert spaces}
A (\textit{complex}) \textit{Hilbert space} is a vector space \(X\) equipped with a complete inner product \(\innerProd{\blank}{\blank} \colon X \times X \to \Comps\). In this article, inner products will be linear in their second argument. Intuitively, inner products encode geometry. Concretely, the length of a vector \(x\) is given by its norm \(\norm{x} = \sqrt{\innerProd{x}{x}}\) and the angle \(\theta\) between two non-zero vectors \(x\) and \(y\) is given by the formula
\[\cos \theta = \frac{\Re \innerProd{x}{y}}{\norm{x}\norm{y}}.\]
Prototypical examples of Hilbert spaces include the space \(\Comps^n\) with the inner product
\[\innerProd[\big]{(x_1, x_2, \dots, x_n)}{(y_1, y_2, \dots, y_n)} = {x_1}^{\!\inv} y_1 + {x_2}^{\!\inv} y_2 + \dots + {x_n}^\inv y_n,\]
where \(r^\inv\) denotes the conjugate of a complex number \(r\), and the space
\[\ell_2(\Nats) = \setb[\big]{(x_1, x_2, \dots) \in \Comps^{\Nats}}{\abs{x_1}^2 + \abs{x_2}^2 + \cdots < \infty}\]
of square-summable sequences with the inner product
\[\innerProd[\big]{(x_1, x_2, \dots)}{(y_1, y_2, \dots)} = {x_1}^{\!\inv} y_1 + {x_2}^{\!\inv} y_2 + \cdots.\]

A function \(f \colon X \to Y\) between Hilbert spaces is called \textit{bounded} if there is a real number \(C \geq 0\) such that \(\norm{fx} \leq C \norm{x}\) for each \(x \in X\). A linear map between Hilbert spaces is bounded if and only if it is continuous with respect to the topologies induced by the norms. Every linear map between finite-dimensional Hilbert spaces is bounded. Write \(\Hilb\) for the category of Hilbert spaces and bounded linear maps.

There is a bijection between the vectors in a Hilbert space \(X\) and the morphisms \(\Comps \to X\) in \(\Hilb\). One direction of the bijection sends a vector \(x \in X\) to the bounded linear map \(r \mapsto xr\). The other direction sends a morphism \(f \colon \Comps \to X\) to the vector~\(f1\). For notational convenience, we will use the same name for vectors in \(X\) and the corresponding morphisms \(\Comps \to X\) in \(\Hilb\). This of course also applies when \(X = \Comps\).

\subsection{Involution}
A \textit{\(\inv\)\nobreakdash-category} is a category equipped with a choice of morphism \(f^\inv \colon Y \to X\) for each morphism \(f \colon X \to Y\), such that
\begin{align*}
    1^\inv &= 1, &
    (gf)^\inv &= f^\inv g^\inv, &
    &\text{and} &
    (f^\inv)^\inv &= f.
\end{align*}
The operation \(f \mapsto f^\inv\) will be called the \textit{involution} of the \(\inv\)\nobreakdash-category. Viewed as a contravariant functor, it exhibits an equivalence between the underlying category and its opposite. We will use this fact repeatedly throughout the article.

The prototypical \(\inv\)\nobreakdash-category is \(\Hilb\). The involution sends a morphism \(f \colon X \to Y\) to its \textit{adjoint}, that is, the unique function \(g \colon Y \to X\) such that \(\innerProd{y}{fx} = \innerProd{gy}{x}\) for all \(x \in X\) and all \(y \in Y\). For each morphism \(f \colon \Comps^m \to \Comps^n\) in \(\Hilb\), with respect to the standard bases, the matrices of \(f\) and \(f^\inv\) are conjugate-transposes of each other.

The inner product of a Hilbert space may be recovered from adjoints via the formula \(\innerProd{x}{y} = x^\inv y\). Just as the inner product of a Hilbert space encodes the geometry of the space, the involution of a \(\inv\)\nobreakdash-category endows the morphisms of the category with some kind of geometry. In particular, the concepts of being \textit{orthogonal} and \textit{isometric}, which are prominent in all aspects of \(\inv\)\nobreakdash-category theory, have a geometric flavour.

Morphisms \(f \colon X \to Z\) and \(g \colon Y \to Z\) are \textit{orthogonal} if \(f^\inv g = 0\). Morphisms \(f \colon X \to Z\) and \(g \colon Y \to Z\) in \(\Hilb\) are orthogonal exactly when their ranges are orthogonal as subspaces of \(Z\), that is, when \(\innerProd{fx}{gy} = 0\) for all \(x \in X\) and \(y \in Y\).

A morphism~\(f\) is \textit{isometric} if \(f^\inv \! f = 1\). A morphism \(f \colon X \to Y\) in \(\Hilb\) is isometric if and only if \(\norm{fx} = \norm{x}\) for all \(x \in X\). An \textit{isometry} is a morphism that is isometric. The terms \textit{coisometric} and \textit{coisometry} are defined dually. A morphism is \textit{unitary} if it is isometric and coisometric, or, equivalently, if it is isometric and invertible.\footnote{In the dagger\nobreakdash-category literature, \textit{dagger monomorphism}, \textit{dagger epimorphism} and \textit{dagger isomorphism} are synonyms for \textit{isometry}, \textit{coisometry}, and \textit{unitary morphism}, respectively.} Identity morphisms are unitary, and the classes of isometries and coisometries are both closed under composition, so all three classes of morphisms form wide subcategories. For each \(\inv\)\nobreakdash-category \(\C\), the wide subcategory of isometries of \(\C\) will be denoted \(\NMono{\C}\).

A morphism \(f\) is \textit{partially isometric} if \(f = ff^\inv \! f\). A \textit{partial isometry} is a morphism that is partially isometric. All isometries and coisometries are partially isometric. If \(f\) is partially isometric, then \(f^\inv\) is also partially isometric. Partially isometric morphisms are not, in general, closed under composition.

If a morphism \(f\) is invertible, then \((f^{-1})^\inv = (f^\inv)^{-1}\). An endomorphism \(f\) is called \textit{Hermitian} if \(f^\inv = f\).

\subsection{Closed subobjects}
\label{sect:closed-subobject}
There is a bijection between the closed subspaces of a Hilbert space \(X\) and the subobjects of \(X\) in \(\Hilb\) that are represented by an isometry. Indeed, if \(A\) is a closed subspace of \(X\), then it is a Hilbert space with respect to the restriction of the inner product of \(X\), and the inclusion map \(A \hookrightarrow X\) is isometric. Conversely, if \(m \colon A \to X\) is an isometry, then its range is a closed subspace of \(X\). Beware: subobjects with an isometric representative may also have other representatives that are \textit{not} isometric. Indeed, isometries are not, in general, preserved by precomposition with isomorphisms, that is, if \(m \colon A \to X\) is an isometry and \(f \colon B \to A\) is an isomorphism, then the monomorphism \(mf\) is not necessarily isometric. The monomorphisms in \(\Hilb\) that represent closed subspaces are actually characterised by the following new notion of closed monomorphism.

\begin{definition}
    \label{def:closed-mono}
    A morphism \(f\) in a \(\inv\)\nobreakdash-category is a \textit{closed monomorphism} if \(f^\inv \! f\) is an isomorphism, and is a \textit{closed epimorphism} if \(ff^\inv\) is an isomorphism. A subobject is \textit{closed} if it has a representative that is a closed monomorphism.
\end{definition}

Generalising from the theory of Hilbert spaces, one might expect that closed subobjects always have isometric representatives. Beware: this is not the case in general. However, it is the case in \(\inv\)\nobreakdash-categories like \(\Hilb\) in which every parallel pair has an isometric equaliser (see \cref{sect:isometric-kernel} below).

\begin{proposition}
\label{prop:closed-mono-retraction}
A morphism \(s \colon A \to X\) is a closed monomorphism if and only if it has a retraction \(r \colon X \to A\) such that the idempotent \(sr\) is Hermitian, in which case \((s^\inv \! s)^{-1} = rr^\inv\) and \(r = (s^\inv \! s)^{-1}s^\inv\).
\end{proposition}

\begin{proof}
If \(s\) is a closed monomorphism, then \((s^\inv \! s)^{-1}s^\inv\) is a retraction of \(s\), and \(s(s^\inv \! s)^{-1}s^\inv\) is Hermitian. Conversely, if \(r \colon X \to A\) is a retraction of \(s\) and \(sr\) is Hermitian, then \(s^\inv \! s\) is an isomorphism with inverse \(rr^\inv\) because
\[s^\inv \! s rr^\inv = s^\inv (sr)^\inv r^\inv = s^\inv r^\inv \! s^\inv r^\inv = (rs)^\inv (rs)^\inv = 1^\inv 1^\inv = 1\]
and
\[rr^\inv \! s^\inv \! s = r(sr)^\inv \! s = rsrs = 1.\]
Finally
\[(s^\inv \! s)^{-1}s^\inv = rr^\inv \! s^\inv = r(sr)^\inv = rsr = 1r = r.\qedhere\]
\end{proof}

We will refer to \(r = (s^\inv \! s)^{-1}s^\inv\) as the \textit{canonical retraction} of \(s\). Beware: there may be other retractions of \(s\). In \(\Hilb\), if \(r\) is the canonical retraction of \(s\) then \(sr\) is the orthogonal projection onto the image of \(s\), while if \(r\) is some other retraction of \(s\) then \(sr\) is some non-orthogonal projection onto the image of \(s\).

By \cref{prop:closed-mono-retraction}, closed monomorphisms are split monomorphisms, and thus are indeed monic. Dually, closed epimorphisms are split epimorphisms, and thus are epic. A morphism \(s\) is a closed monomorphism if and only if \(s^\inv\) is a closed epimorphism. All isometries are closed monomorphisms, all coisometries are closed epimorphisms, and all isomorphisms are both closed monomorphisms and closed epimorphisms.

\begin{remark}
\label{p:closed-mono-iso}
The composite \(sf\) of a closed monomorphism \(s \colon A \to X\) with an isomorphism \(f \colon B \to A\) is again a closed monomorphism. Indeed
\[(f^\inv \! s^\inv \! s f)^{-1} = f^{-1}(s^\inv \! s)^{-1} f^{-1\inv}.\] In particular, every representative of a closed subobject is a closed monomorphism.
\end{remark}

\subsection{Zero objects}
A \textit{zero object} is an object that is both initial and terminal. The symbol \(\zero\) denotes a chosen zero object. In a category with a zero object, for all objects \(X\) and \(Y\), there is a unique morphism from \(X\) to \(Y\) that factors through a zero object; it is called the \textit{zero morphism} from \(X\) to \(Y\) and denoted by \(0 \colon X \to Y\). Zero morphisms satisfy the annihilation laws \(0f = 0\) and \(f0 = 0\). There is only one morphism \(\zero \to \zero\); in particular, the morphisms \(0 \colon \zero \to \zero\) and \(1 \colon \zero \to \zero\) coincide. 

In a \(\inv\)\nobreakdash-category, the concepts of initial object, terminal object and zero object coincide. Also, zero morphisms satisfy the equation \(0^\inv = 0\). Additionally, for each object \(X\), the morphism \(0 \colon \zero \to X\) is isometric. The category \(\Hilb\) has a canonical zero object, namely \(\zero = \set{0}\). For all objects \(X\) and~\(Y\) of \(\Hilb\), the zero morphism \(0 \colon X \to Y\) is the function that sends every element of~\(X\) to the vector \(0\) in \(Y\).

\subsection{Isometric kernels}
\label{sect:isometric-kernel}
In a category with a zero object, a \textit{kernel} of a morphism \(f \colon X \to Y\) is an equaliser of \(f\) and \(0 \colon X \to Y\). If such a kernel exists, the subobject that it represents will be denoted \(\Ker f\). A \textit{normal monomorphism} is a morphism that is a kernel of some morphism. Every normal monomorphism is indeed a monomorphism. The terms \textit{cokernel} and \textit{normal epimorphism} are defined dually.

In a \(\inv\)\nobreakdash-category, a morphism \(m\) is an equaliser of \(f\) and \(g\) if and only if \(m^\inv\) is a coequaliser of \(f^\inv\) and \(g^\inv\). Hence, in a \(\inv\)\nobreakdash-category, every parallel pair has an equaliser if and only if every parallel pair has a coequaliser. Similarly, in a \(\inv\)\nobreakdash-category with a zero object, every morphism has a kernel if and only if every morphism has a cokernel.

In \(\inv\)\nobreakdash-category theory, kernels and equalisers that are closed monomorphisms or isometries are of particular interest.\footnote{In the dagger\nobreakdash-category literature, isometric equalisers and isometric kernels are referred to as \textit{dagger equalisers} and \textit{dagger kernels}. Dually, coisometric coequalisers and coisometric cokernels are referred to as \textit{dagger coequalisers} and \textit{dagger cokernels}.} In a \(\inv\)\nobreakdash-category with a zero object, if a morphism \(f\) has an isometric kernel, then, by \cref{p:closed-mono-iso}, every kernel of \(f\) is a closed monomorphism (but not necessarily an isometry). In particular, if every morphism has an isometric kernel, then every normal monomorphism is a closed monomorphism. We revisit these observations in \cref{lem:split-is-closed}.

In \(\Hilb\), all parallel pairs of morphisms have an isometric equaliser, and thus all morphisms have an isometric kernel. If \(f \colon X \to Y\) is a bounded linear map, then 
\[\Ker f = \setb[\big]{x \in X}{fx = 0}\]
is a closed subspace of \(X\), so the inclusion map \(\Ker f \hookrightarrow X\) is isometric, and is actually an isometric kernel of \(f\) in \(\Hilb\). If \(g \colon X \to Y\) is another bounded linear map, then the inclusion into \(X\) of the closed subspace
\[\Eq(f, g) = \setb[\big]{x \in X}{fx = gx}\]
of \(X\) is similarly an isometric equaliser of \(f\) and \(g\) in \(\Hilb\).

\subsection{Orthonormal biproducts}
\label{sect:orthonormal-biproducts}
In a category with a zero object, a \textit{biproduct} of a finite family \(X_1, \dots, X_n\) of objects is a tuple \((X, s_1, r_1, \dots, s_n, r_n)\) consisting of an object \(X\) and morphisms \(s_j \colon X_j \to X\) and \(r_j \colon X \to X_j\) for each \(j \in \set{1, \dots, n}\), such that \((X, s_1, \dots, s_n)\) is a coproduct, \((X, r_1, \dots, r_n)\) is a product, and
\begin{equation}
    \label{eq:biproduct-retract}
    r_k s_j = \begin{cases} 1 &\text{if } j = k,\\0 &\text{otherwise,} \end{cases}
\end{equation}
for all \(j,k \in \set{1, \dots, n}\). A biproduct of an empty family of objects is precisely a zero object. A category \textit{has finite biproducts} if it has a zero object and every finite family of objects has a biproduct, or, equivalently, if it has a zero object and every pair of objects has a biproduct.
The tuple \((X_1 \oplus \dots \oplus X_n,\, i_1, p_1, \dots, i_n, p_n)\) is a chosen biproduct of \(X_1, \dots, X_n\). The chosen biproduct of \(n\) copies of \(X\) is written \(X^{\oplus n}\).

Given morphisms
\(f_{kj} \colon X_j \to Y_k\)
for each \(j \in \set{1, \dots, m}\) and each \(k \in \set{1, \dots, n}\), there is a unique morphism
\[
    \begin{bmatrix}
        f_{11} & \cdots & f_{1m}\\
        \vdots & \ddots & \vdots\\
        f_{n1} & \cdots & f_{nm}
    \end{bmatrix}
    \colon X_1 \oplus \dots \oplus X_m \to Y_1 \oplus \dots \oplus Y_n
\]
such that
\[
    p_k
    \begin{bmatrix}
        f_{11} & \cdots & f_{1m}\\
        \vdots & \ddots & \vdots\\
        f_{n1} & \cdots & f_{nm}
    \end{bmatrix}
    i_j
    = f_{kj}
\]
for each \(j \in \set{1, \dots, m}\) and each \(k \in \set{1, \dots, n}\).

The \textit{diagonal} \(\diagonal \colon X \to X^{\oplus n}\) and the \textit{codiagonal} \(\codiagonal \colon X^{\oplus n} \to X\) are defined by
\[
    \diagonal = \begin{bmatrix}1 \\ \vdots \\ 1\end{bmatrix}
    \qquad\text{and}\qquad
    \codiagonal = \begin{bmatrix}1 & \cdots & 1\end{bmatrix}.
\]

A choice of finite biproducts in a category \(\C\) extends to functors \[\underbrace{- \oplus - \oplus \dots \oplus -}_{n \text{ times}} \colon \underbrace{\C \times \C \times \dots \times \C}_{n \text{ times}} \to \C\]
for each natural number \(n\), whose actions on morphisms are defined by the equation
\[
    f_1 \oplus f_2 \oplus \dots \oplus f_n
    = \begin{bmatrix}
        f_1 & 0 & \cdots & 0\\
        0 & f_2 & \cdots & 0\\
        \vdots & \vdots & \ddots & \vdots\\
        0 & 0 & \cdots & f_n
    \end{bmatrix}.
\]

A coproduct \((X, s_1, \dots, s_n)\) in a \(\inv\)\nobreakdash-category is \textit{orthonormal} if \(s_1, \dots, s_n\) are pairwise-orthogonal isometries. A biproduct \((X, s_1, r_1, \dots, s_n, r_n)\) in a \(\inv\)\nobreakdash-category is \textit{orthonormal} if \(r_k = {s_k}^{\!\inv}\) for each \(k\).\footnote{The dagger\nobreakdash-category literature uses the names \textit{dagger coproduct} and \textit{dagger biproduct}.} If the tuple \((X, s_1, \dots, s_n)\) is an orthonormal coproduct, then \((X, s_1, {s_1}^\inv, \dots, s_n, {s_n}^\inv)\) is an orthonormal biproduct. Conversely, if \((X, s_1, r_1, \dots, s_n, r_n)\) is an orthonormal biproduct, then \((X, s_1, \dots, s_n)\) is an orthonormal coproduct. Orthonormal coproduct and orthonormal biproduct are thus different perspectives on a single concept. \textit{Orthogonal biproducts}---a weaker variant of this concept---are introduced in \cref{sect:gram-schmidt} for the generalised Gram--Schmidt process. The adjectives \textit{orthogonal} and \textit{orthonormal} are of course lifted from linear algebra. Given an orthogonal (or orthonormal) basis of a finite-dimensional Hilbert space \(X\), the bounded linear maps \(\Comps \to X\) associated to each vector in the basis form an orthogonal (or orthonormal) coproduct in \(\Hilb\).

A \(\inv\)\nobreakdash-category \textit{has finite orthonormal biproducts} if it has a zero object and every finite family of objects has an orthonormal biproduct, or, equivalently, if it has a zero object and every pair of objects has an orthonormal biproduct. In a \(\inv\)\nobreakdash-category that has finite orthonormal biproducts, the chosen biproducts \(X_1 \oplus \dots \oplus X_n\) are assumed to be orthonormal.

The \(\inv\)\nobreakdash-category \(\Hilb\) has finite orthonormal biproducts. An orthonormal coproduct of Hilbert spaces \(X_1, \dots, X_n\) may be constructed from their direct sum \((X_1 \oplus \dots \oplus X_n,\, i_1, \dots, i_n)\) as vector spaces by giving \(X_1 \oplus \dots \oplus X_n\) the inner product 
\[\innerProd[\big]{(x_1, \dots, x_n)}{(y_1, \dots, y_n)} = \innerProd{x_1}{y_1} + \dots + \innerProd{x_n}{y_n}.\]

Given morphisms
\(f_{kj} \colon X_j \to Y_k\)
for each \(j \in \set{1, \dots, m}\) and each \(k \in \set{1, \dots, n}\), if the biproducts \(X_1 \oplus \dots \oplus X_m\) and \(Y_1 \oplus \dots \oplus Y_n\) are orthonormal, then
\begin{equation}
    \label{eq:adjoint-matrix-orthonormal}
    \begin{bmatrix}
        f_{11} & \cdots & f_{1m}\\
        \vdots & \ddots & \vdots\\
        f_{n1} & \cdots & f_{nm}
    \end{bmatrix}^\inv
    =
    \begin{bmatrix}
        {f_{11}}^\inv & \cdots & {f_{n1}}^\inv\\
        \vdots & \ddots & \vdots\\
        {f_{1m}}^\inv & \cdots & {f_{nm}}^\inv
    \end{bmatrix}.
\end{equation}
It follows also that
\[\diagonal^\inv = \codiagonal \qquad\text{and}\qquad (f_1 \oplus f_2 \oplus \dots \oplus f_n)^\inv = {f_1}^{\!\inv} \oplus {f_2}^{\!\inv} \oplus \dots \oplus {f_n}^\inv.\]

\subsection{Pre-Hilbert \texorpdfstring{\(\inv\)}{*}-categories}
\label{subsect:r-cat}
We have now encountered all of the concepts that we need to understand the definition of pre\nobreakdash-Hilbert \(\inv\)\nobreakdash-category, which is recalled here from the introduction.

\begin{definition}
    \label{def:r-category}
    A \textit{pre\nobreakdash-Hilbert \(\inv\)\nobreakdash-category} is a \(\inv\)\nobreakdash-category in which
    \begin{enumerate}[(R1)]
        \item 
        \label{axiom:R-zero}
        there is a zero object,
        \item
        \label{axiom:R-coproduct}
        every pair of objects has an orthonormal biproduct,
        \item
        \label{axiom:R-kernel}
        every morphism has an isometric kernel,
        \item
        \label{axiom:R-normal}
        every diagonal \(\diagonal \colon X \to X \oplus X\) is a normal monomorphism.
    \end{enumerate}
\end{definition}

We have also seen that the \(\inv\)\nobreakdash-category \(\Hilb_{\Comps} = \Hilb\) of complex Hilbert spaces and bounded linear maps satisfies axioms \cref{axiom:R-zero,axiom:R-kernel,axiom:R-coproduct}. For axiom \cref{axiom:R-normal}, the diagonal \(\diagonal \colon X \to X \oplus X\) is a kernel of the morphism \(X \oplus X \to X\) defined by \((x, y) \mapsto x - y\). The analogous \(\inv\)\nobreakdash-categories \(\Hilb_{\Reals}\) and \(\Hilb_{\Quats}\) of real and quaternionic Hilbert spaces are also pre\nobreakdash-Hilbert \(\inv\)\nobreakdash-categories by similar reasoning. More examples of pre\nobreakdash-Hilbert \(\inv\)\nobreakdash-categories of an inner product flavour are described in detail in \cref{sect:examples}.

Given a small category \(\C\) and a (locally small) category \(\D\), write \([\C, \D]\) for the (locally small) category whose objects are the functors from \(\C\) to \(\D\), and whose morphisms are the natural transformations between them. If \(\D\) is an abelian category, then so is the functor category \([\C, \D]\). In the rest of this section, we prove a similar fact about pre\nobreakdash-Hilbert \(\inv\)\nobreakdash-categories, building on the work of \textcite{karvonen:limits-dagger-cats}. 

A \textit{\(\inv\)\nobreakdash-functor} is a functor \(F\) between \(\inv\)\nobreakdash-categories such that \((F f)^\inv = F f^\inv\) for all morphisms \(f\) in the domain of \(F\). If \(F\) and \(G\) are \(\inv\)\nobreakdash-functors, then applying the involution componentwise to a natural transformation \(\sigma \colon F \Rightarrow G\) yields a natural transformation \(\sigma^\inv \colon G \Rightarrow F\). Given a small \(\inv\)\nobreakdash-category \(\C\) and a \(\inv\)\nobreakdash-category \(\D\), write \([\C, \D]_\inv\) for the full subcategory of \([\C, \D]\) spanned by the \(\inv\)\nobreakdash-functors. It is a \(\inv\)\nobreakdash-category when equipped with the pointwise involution described above~\cite[Example~2.7]{karvonen:limits-dagger-cats}.

\begin{proposition}
For each small \(\inv\)\nobreakdash-category \(\C\) and each pre\nobreakdash-Hilbert \(\inv\)\nobreakdash-category \(\D\), the \(\inv\)\nobreakdash-category \([\C, \D]_\inv\) is also a pre\nobreakdash-Hilbert \(\inv\)\nobreakdash-category.
\end{proposition}

\begin{proof}
Similar to the way that limits in \([\C, \D]\) are computed componentwise~\cite[Proposition~3.3.9]{riehl:category-theory-context}, orthonormal biproducts and isometric kernels in \([\C, \D]_\inv\) are also computed componentwise~\cite[Examples~3.2 and 3.8]{karvonen:limits-dagger-cats}. Hence, as \(\D\) satisfies axioms \cref{axiom:R-zero,axiom:R-coproduct,axiom:R-kernel}, so does \([\C, \D]_{\inv}\). Proving that \([\C, \D]_\inv\) also satisfies axiom \cref{axiom:R-normal} requires more work. Let \(F \colon \C \to \D\) be a \(\inv\)\nobreakdash-functor. As orthonormal biproducts in \([\C, \D]_\inv\) are computed componentwise, for each object \(X\) of \(\C\), the \(X\)-component of the diagonal \(\Delta \colon F \Rightarrow F \oplus F\) in \([\C, \D]_\inv\) is the diagonal \(\Delta_X \colon FX \to FX \oplus FX\) in~\(\D\). By the dual of axiom \cref{axiom:R-kernel}, the morphism \(\Delta_X\) has a coisometric cokernel \(\sigma_X \colon FX \oplus FX \to GX\). Universality of these cokernels yields an extension of this data to a functor \(G \colon \C \to \D\) and a natural transformation \(\sigma \colon F \oplus F \Rightarrow G\). Now \(F \oplus F\) is a \(\inv\)\nobreakdash-functor and the components of \(\sigma\) are coisometries, so \(G\) is actually also a \(\inv\)\nobreakdash-functor~\cite[Lemma~2.12]{karvonen:limits-dagger-cats}. Now \(\Delta_X\) is a kernel of \(\sigma_X\) in \(\D\) for all objects \(X\) of \(\C\), and limits in \([\C, \D]\) are computed componentwise, so \(\Delta\) is actually a kernel of \(\sigma\) in \([\C, \D]\). But \([\C, \D]_\inv\) is a full subcategory of \([\C, \D]\) and the limit object \(F\) is a \(\inv\)\nobreakdash-functor, so \(\Delta\) is also a kernel of \(\sigma\) in \([\C, \D]_\inv\)~\cite[Lemma~3.3.5]{riehl:category-theory-context}. Hence \([\C, \D]_\inv\) satisfies axiom~\cref{axiom:R-normal}.
\end{proof}

A \textit{groupoid} is a small category in which every morphism is invertible. A \textit{group} is merely a groupoid with one object. Groupoids are canonically \(\inv\)\nobreakdash-categories with \(f^\inv = f^{-1}\) for each morphism \(f\). A \textit{unitary representation} of a groupoid \(\C\) is a \(\inv\)\nobreakdash-functor \(F \colon \C \to \Hilb_{\Comps}\)~\cite[Example~2.7]{karvonen:limits-dagger-cats}. Concretely, it consists of a choice of Hilbert spaces \(FX\) for all objects \(X\) of \(\C\) and a functorial choice of unitary morphisms \(Ff \colon FX \to FY\) for all morphisms \(f \colon X \to Y\) of \(\C\). Given unitary representations \(F\) and \(G\), an \textit{intertwiner} from \(F\) to \(G\) is a natural transformation \(\sigma \colon F \Rightarrow G\).

\begin{corollary}
The category of unitary representations of a groupoid and their intertwiners is a pre\nobreakdash-Hilbert \(\inv\)\nobreakdash-category.
\end{corollary}

\section{Similarities with abelian categories}
\label{sect:abelian-properties}

This section is about the properties of pre\nobreakdash-Hilbert \(\inv\)\nobreakdash-categories that resemble well-known properties of abelian categories. Some of these follow merely from additivity---the focus of \cref{sec:additivity}. There, we will also see two alternative definitions of pre\nobreakdash-Hilbert \(\inv\)\nobreakdash-category that are, by additivity, equivalent to the original one. Many of the other abelian-like properties of pre\nobreakdash-Hilbert \(\inv\)\nobreakdash-categories follow from \textit{quasi-abelianness}---the focus of \cref{sec:quasiabelian}.

\subsection{Enrichment in commutative monoids}
A \textit{category enriched in commutative monoids} is a category whose hom-sets each come equipped with a binary operation \(+\) and a distinguished element \(0\) that make them into commutative monoids, such that
\[
    f0 = 0,
    \qquad
    0g = 0,
    \qquad
    f(g_1 + g_2) = fg_1 + fg_2,
    \qquad
    (f_1 + f_2)g = f_1g + f_2g.
\]
In such a category, a tuple \((X, s_1, r_1, \dots, s_n, r_n)\) is a biproduct of \((X_1, \dots, X_n)\) if and only if equation \cref{eq:biproduct-retract} holds for all \(j, k \in \set{1, \dots, n}\) and
\begin{equation}
    \label{eq:biproduct-span}
    s_1 r_1 + \dots + s_n r_n = 1.
\end{equation} 
This is because equation \cref{eq:biproduct-span} ensures that the tuple \((X, s_1, \dots, s_n)\) is a coproduct and that the tuple \((X, r_1, \dots, r_n)\) is a product.

If a category has finite biproducts then it also has a unique enrichment in commutative monoids~\cite[Lemma~2.4]{vicary:2011:complex-numbers}. The distinguished element \(0 \colon X \to Y\) is the zero morphism from \(X\) to \(Y\), and
\begin{equation}
    \label{eq:sum-formula}
    f_1 + \dots + f_n = \codiagonal \begin{bmatrix}f_1 \\ \vdots \\ f_n \end{bmatrix} = \codiagonal(f_1 \oplus \dots \oplus f_n) \diagonal =  \begin{bmatrix}f_1 & \cdots & f_n \end{bmatrix} \diagonal.
\end{equation}
With respect to this enrichment in commutative monoids, the usual formula 
\begin{equation}
\label{eq:matrix-product}
    p_\ell
    \begin{bmatrix}
        g_{11} & \cdots & g_{1n}\\
        \vdots & \ddots & \vdots\\
        g_{p1} & \cdots & g_{pn}
    \end{bmatrix}
    \begin{bmatrix}
        f_{11} & \cdots & f_{1m}\\
        \vdots & \ddots & \vdots\\
        f_{n1} & \cdots & f_{nm}
    \end{bmatrix}
    i_j = g_{\ell 1}f_{1j} + \dots + g_{\ell n} f_{nj}
\end{equation}
for the product of two matrices holds, whenever it is well defined.

It is well known that in \(\inv\)\nobreakdash-categories that have finite \textit{orthonormal} biproducts, the involution respects addition. Surprisingly, orthonormality is not necessary.

\begin{proposition}
\label{prop:star-addition}
In a \(\inv\)\nobreakdash-category that has finite biproducts, for all \(f, g \colon X \to Y\),
\[(f + g)^\inv = f^\inv + g^\inv.\]
\end{proposition}

\begin{proof}
Observe that \(\diagonal^\inv \colon X \oplus X \to X\) is the codiagonal on \(X\) with respect to the biproduct \((X \oplus X,\, {p_1}^{\!\inv}\!,\, {i_1}^{\!\inv}\!,\, {p_2}^{\!\inv}\!, \,{i_2}^{\!\inv})\), and \(\copair{f}{g}^\inv \colon Y \to X \oplus X\) is the product pairing of \(f^\inv\) and \(g^\inv\) with respect to this biproduct. Hence, applying equation \cref{eq:sum-formula} twice, first with respect to the biproduct \((X \oplus X, \,i_1, p_1, i_2, p_2)\) and then with respect to the biproduct \((X \oplus X,\, {p_1}^{\!\inv}\!,\, {i_1}^{\!\inv}\!,\, {p_2}^{\!\inv}\!, \,{i_2}^{\!\inv})\), we see that
\[(f + g)^\inv = \paren[\big]{\copairBig{f}{g} \diagonal}^\inv = \diagonal^\inv \copairBig{f}{g}^\inv = f^\inv + g^\inv. \qedhere\] 
\end{proof}

\subsection{Additivity}
\label{sec:additivity}
A \textit{category enriched in abelian groups} is a category enriched in commutative monoids in which the commutative monoid structure on each hom-set is actually an abelian group. In a category enriched in commutative monoids, an object \(X\) is \textit{abelian} if \(1 \colon X \to X\) has an additive inverse~(see \cite[Proposition~1.5.3]{borceux-bourn:2004:malcev-protomodular-homological}). A category enriched in commutative monoids is enriched in abelian groups if and only if every object is abelian. A category is \textit{additive} if it has finite biproducts and its unique enrichment in commutative monoids is actually an enrichment in abelian groups, that is, if it has finite biproducts and every object is abelian.

\begin{proposition}
\label{cor:additive}
Every pre\nobreakdash-Hilbert \(\inv\)\nobreakdash-category is additive.
\end{proposition}

We will prove \cref{cor:additive} using the following new characterisation of abelian objects in categories that have finite biproducts.

\begin{theorem}
\label{prop:semiadditive-abelian-objects}
In a category that has finite biproducts, an object \(X\) is abelian if and only if the diagonal \(\diagonal \colon X \to X \oplus X\) is a kernel of a split epimorphism.
\end{theorem}

\begin{proof}
Suppose that \(X\) is abelian. Since the morphism \(\copair{1}{-1} \colon X \oplus X \to X\) has a section, namely \(i_1 \colon X \to X \oplus X\), it suffices to show that \(\diagonal \colon X \to X \oplus X\) is a kernel of \(\copair{1}{-1}\). First, \(\copair{1}{-1}\diagonal = 1 + (-1) = 0\). Next, consider a morphism \(f \colon Y \to X \oplus X\) such that \(\copair{1}{-1}f = 0\). Then \(p_1 f = p_2 f\) because \(p_1f - p_2 f = \copair{1}{-1}f = 0\). Hence
\[
    f
    = i_1p_1f + i_2p_2f
    = i_1p_1f + i_2p_1f
    = (i_1 + i_2)p_1f
    = \diagonal p_1f.
\]
This factorisation is unique because \(p_1\diagonal = 1\).

Conversely, suppose that \(\diagonal \colon X \to X \oplus X\) is a kernel of a morphism \(r \colon X \oplus X \to A\) with a section \(s \colon A \to X \oplus X\), and let \(e = sr\) be the induced idempotent.

The fact that \(\diagonal\) is a kernel of \(r\) may be reinterpreted in terms of \(e\) as follows:
\[
    e(i_1 + i_2) = 0,
    \qquad\text{and}\qquad
    ef = 0 \implies p_1f = p_2f.
\]
Indeed, \(e(i_1 + i_2) = e\diagonal = sr\diagonal = s0 = 0\). Also, if \(f \colon Y \to X \oplus X\) satisfies \(ef = 0\), then
\(rf = rsrf = ref = r0 = 0\); by universality of the kernel, there is a unique morphism \(f' \colon Y \to X\) such that \(f = \diagonal f'\); hence \(p_1f = p_1\diagonal f' = 1f' = p_2\diagonal f' = p_2 f\).

We will show that \(p_1ei_2 + p_2ei_1\) is an additive inverse of \(1 \colon X \to X\). Observe that
\begin{align*}
    ei_1(1 + p_1ei_2 + p_2ei_1)
    &= ei_1 + ei_1p_1ei_2 + ei_1p_2ei_1
    \\&= ei_1 + ei_1p_1ei_2 + ei_1p_2ei_1 + 0p_2ei_2
    \\&= ei_1 + ei_1p_1ei_2 + ei_1p_2ei_1 + e(i_1 + i_2)p_2ei_2
    \\&= ei_1 + e(i_1p_1 + i_2p_2)ei_2 + ei_1p_2e(i_1 + i_2)
    \\&= ei_1 + e1ei_2 + ei_1p_20
    \\&= ei_1 + e^2i_2
    \\&= e(i_1 + i_2)
    \\&= 0.
    \\ \intertext{From the implication proved earlier, it follows that}
    1 + p_1ei_2 + p_2ei_1
    &= p_1i_1(1 + p_1ei_2 + p_2ei_1)
    \\&= p_2i_1(1 + p_1ei_2 + p_2ei_1)
    \\&= 0(1 + p_1ei_2 + p_2ei_1)
    \\&= 0. \qedhere
\end{align*}
\end{proof}

\begin{proof}[Proof of \cref{cor:additive}]
As all pre\nobreakdash-Hilbert \(\inv\)\nobreakdash-categories have finite biproducts, it suffices to show that every object of a pre\nobreakdash-Hilbert \(\inv\)\nobreakdash-category is abelian. Let \(X\) be such an object. Then the diagonal \(\diagonal \colon X \to X \oplus X\) is a normal monomorphism by axiom \cref{axiom:R-normal}, and it has a coisometric cokernel by the dual of axiom~\cref{axiom:R-kernel}. As the diagonal \(\diagonal \colon X \to X \oplus X\) is normal monomorphism, it is a kernel of its cokernel. By \cref{prop:semiadditive-abelian-objects}, the object \(X\) is abelian.
\end{proof}

Additivity can, instead of being a proposition, be incorporated into the definition.

\begin{proposition}
\label{prop:definition-by-abelian}
A \(\inv\)\nobreakdash-category is a pre\nobreakdash-Hilbert \(\inv\)\nobreakdash-category if and only if it has an enrichment in abelian groups and satisfies axioms \cref{axiom:R-zero}, \cref{axiom:R-coproduct} and \cref{axiom:R-kernel}.
\end{proposition}

\begin{proof}
For the \textit{if} direction, the diagonal \(\diagonal \colon X \to X \oplus X\) is a kernel of the morphism \(\copair{1}{-1} \colon X \oplus X \to X\). The \textit{only if} direction is precisely \cref{cor:additive}.
\end{proof}

While the original definition is nicer aesthetically---it refers only to a few simple properties of products and kernels---this alternative definition is easier to work with. It is used in \cref{sect:examples} to show that certain \(\inv\)\nobreakdash-categories are pre\nobreakdash-Hilbert \(\inv\)\nobreakdash-categories.

As promised in the introduction, pre\nobreakdash-Hilbert \(\inv\)\nobreakdash-categories are also equivalently defined by a subset of the axioms in the characterisation~\cite{heunen:2022:axioms-for-the-category} of \(\Hilb_{\field{K}}\).

\begin{proposition}
\label{prop:equivalent-definition-heunen}
A \(\inv\)\nobreakdash-category is a pre\nobreakdash-Hilbert \(\inv\)\nobreakdash-category if and only if
\begin{enumerate}[(H1)]
    \item
    \label{axiom:hilbert-zero}
    there is a zero object,
    \item
    \label{axiom:hilbert-product}
    every pair of objects has an orthonormal biproduct,
    \item
    \label{axiom:hilbert-equaliser}
    every parallel pair of morphisms has an isometric equaliser, and
    \item
    \label{axiom:hilbert-normal}
    every isometry is a normal monomorphism.
\end{enumerate}
\end{proposition}

\begin{proof}
Axioms \cref{axiom:hilbert-zero,axiom:hilbert-product} are identical to axioms \cref{axiom:R-zero,axiom:R-coproduct}.

For the \textit{if} direction, axiom \cref{axiom:R-kernel} is a special case of axiom \cref{axiom:hilbert-equaliser}. For axiom \cref{axiom:R-normal}, the morphisms \(p_1, p_2 \colon X \oplus X \to X\) have an isometric equaliser \(m \colon A \to X \oplus X\) by axiom \cref{axiom:hilbert-equaliser}, and \(m\) is the kernel of a morphism \(f \colon X \oplus X \to Y\) by axiom~\cref{axiom:hilbert-normal}. The diagonal \(\diagonal \colon X \to X \oplus X\) is another equaliser of \(p_1\) and \(p_2\), so it represents the same subobject as \(m\), and thus must be another kernel of \(f\). 

For the \textit{only if} direction, pre\nobreakdash-Hilbert \(\inv\)\nobreakdash-categories are additive by \cref{cor:additive}. Axiom \cref{axiom:hilbert-equaliser} follows from axiom \cref{axiom:R-kernel} because a morphism \(m\) is an equaliser of morphisms \(f\) and \(g\) if and only if \(m\) is a kernel of \(f - g\). Axiom \cref{axiom:hilbert-normal} holds because every isometry \(m\) is a kernel of \(1 - mm^\inv\).
\end{proof}

\subsection{Quasi-abelianness}
\label{sec:quasiabelian}
A category is \textit{quasi-abelian}~\cite{schneiders:1999:quasi-abelian-sheaves} (or \textit{almost abelian}~\cite{rump:2001:almost-abelian}) if it is additive, all morphisms have a kernel and a cokernel, and normal monomorphisms and normal epimorphisms are stable under pushout and pullback, respectively.

\begin{proposition}
\label{prop:quasi-abelian}
Every pre\nobreakdash-Hilbert \(\inv\)\nobreakdash-category is quasi-abelian.
\end{proposition}

\begin{lemma}
\label{lem:split-is-closed}
In a pre\nobreakdash-Hilbert \(\inv\)\nobreakdash-category, the classes of normal monomorphisms, closed monomorphisms and split monomorphisms coincide.
\end{lemma}

\begin{proof}
By \cref{prop:closed-mono-retraction}, every closed monomorphism has a canonical retraction. Also, if \(s\) has a retraction \(r\), then it is a kernel of \(1 - sr\), which exists by \cref{cor:additive}. Finally, if \(s\) is a kernel of a morphism \(f\), then \(f\) has an isometric kernel \(m\) by axiom \cref{axiom:R-kernel}, and \(s = mg\) for some isomorphism \(g\) because kernels are unique up to isomorphism, so \(s\) is a closed monomorphism by \cref{p:closed-mono-iso}.
\end{proof}

We also need the following folklore fact about split monomorphisms.

\begin{lemma}
\label{prop:kernels-pushout-stable}
Split monomorphisms are stable under pushout.
\end{lemma}

\begin{proof}
Let \(s'\) be a pushout of a morphism \(s\) along a morphism \(a\), as depicted in \cref{f:normal-mono-pushout}.
\begin{diagram}
    \centering
    \begin{tikzcd}
        A
            \arrow[r, "s"]
            \arrow[d, "a" swap]
        \&
        X
            \arrow[d]
            \arrow[rdd, "ar", out=-30,in=90]
        \&[-0.9em]
        \\
        A'
            \arrow[r, "{s'}" swap]
            \arrow[drr, "1" swap, out=-45, in=180]
            \&
        |[lrcorner]| X'
            \arrow[dr, "{r'}" near start]
        \&
        \\[-0.6em]
        \&
        \&
        A'
    \end{tikzcd}
    \caption{}
    \label{f:normal-mono-pushout}
\end{diagram}
If \(s\) has a retraction \(r\), then, since \(a = ars\), there is a morphism \(r'\) such that the diagram commutes; in particular, \(r'\) is a retraction of \(s'\).
\end{proof}

\begin{proof}[Proof of \cref{prop:quasi-abelian}]
Additivity is \cref{cor:additive}. Kernels exist by axiom \cref{axiom:R-kernel}. Combining \cref{lem:split-is-closed,prop:kernels-pushout-stable} yields the stability of normal monomorphisms under pushout. The remaining quasi-abelian axioms follow from self-duality.
\end{proof}

Quasi-abelianness is closely connected to two other well-studied generalisations of abelian categories: \textit{regular categories} and \textit{homological categories}. Indeed, the following conditions on a category are equivalent~(see \cite[Remarks~5.3]{rosicky-tholen:2007:factorisation-fibration}):
\begin{enumerate}
    \item it is quasi-abelian;
    \item it is regular, coregular and additive; and
    \item it is homological and cohomological.
\end{enumerate}
Broadly, regular categories recapture some of the exactness properties (compatibility between finite limits and specified types of colimits) of abelian categories without necessarily being additive~\cite{borceux:1994:handbook-categorical-algebra-2}. Homological categories additionally incorporate the classical diagram lemmas, thus abstracting (non-abelian) homological algebra~\cite{borceux-bourn:2004:malcev-protomodular-homological}. In light of \cref{prop:quasi-abelian}, the extensive literature on regular and homological categories is now also a repository of results about pre\nobreakdash-Hilbert \(\inv\)\nobreakdash-categories.

\subsection{Range factorisation}
\label{s:range}

The \textit{range}\footnote{We use the term \textit{range} instead of \textit{image} to avoid a future clash with the notation \(\Im f\) for the imaginary part of a morphism \(f\), and also to avoid confusion with the alternative definition~\cite[Definition~1.1.1]{schneiders:1999:quasi-abelian-sheaves} of the image of \(f\) as the subobject represented by a kernel of a cokernel of \(f\), which is not equivalent to the usual definition of the image of \(f\).} of a morphism \(f \colon X \to Y\), if it exists, is the smallest subobject \(\Ran f\) of \(Y\) through which \(f\) factors~\cite[Section~A1.3]{johnstone:2002:elephant}. A \textit{range factorisation} of a morphism \(f\) is a factorisation of \(f\) through a representative of its range. Every morphism \(f\) in a quasi-abelian category has~\cite[Proposition~1.1.4]{schneiders:1999:quasi-abelian-sheaves} a range factorisation, namely, the unique factorisation of \(f\) through a cokernel of a kernel of~\(f\). In fact, every factorisation \(f = me\) where \(m\) is monic and \(e\) is normal epic is a range factorisation of \(f\); indeed \(f\) and \(e\) have the same kernel because \(m\) is monic, and \(e\) is a cokernel of its kernel because it is normal epic. Since pre\nobreakdash-Hilbert \(\inv\)\nobreakdash-categories are quasi-abelian and have coisometric cokernels, every morphism \(f\) in a pre\nobreakdash-Hilbert \(\inv\)\nobreakdash-category has a range factorisation \(f = me\) where \(e\) is coisometric.

The categorical notion of range does indeed generalise the familiar notion of \textit{range} of a bounded linear map, albeit with some subtleties. The \textit{range} of a bounded linear map \(f \colon X \to Y\) between Hilbert spaces is the subspace
\[\Ran f = \setb[\big]{fx}{x \in X}\]
of \(Y\). With the inner product that it inherits from \(Y\), it is not necessarily even a Hilbert space. A better inner product for \(\Ran f\) is the one transferred from \(X\), which we now construct. Consider the function \((\Ker f)^\perp \to \Ran f\) obtained from \(f\) by restriction and corestriction. It is injective. Indeed, if \(x, x' \in (\Ker f)^\perp\) and \(fx = fx'\), then \(x - x' \in (\Ker f)^\perp \cap \Ker f = \set{0}\), so \(x = x'\). It is also surjective. Indeed, if \(y \in \Ran f\), then \(y = fx\) for some \(x \in X\); now \(x = u + v\) for some \(u \in \Ker f\) and \(v \in (\Ker f)^\perp\), so \(y = fu + fv = fv\). It is thus bijective. By transferring the inner product that \((\Ker f)^\perp\) inherits from \(X\) along this bijection, the vector space \(\Ran f\) becomes a Hilbert space, the corestriction of \(f\) to \(\Ran f\) becomes coisometric, and the inclusion map \(\Ran f \hookrightarrow Y\) becomes bounded.

\section{Orthogonal biproducts and the Gram--Schmidt process}
\label{sect:gram-schmidt}

In a pre\nobreakdash-Hilbert \(\inv\)\nobreakdash-category, if \(X\) is the apex of a coproduct of objects \(X_1, \dots, X_n\), is \(X\) also the apex of an orthonormal coproduct of \(X_1, \dots, X_n\)? Given a finite basis of a Hilbert space, the Gram--Schmidt process produces an orthogonal basis, which can then be normalised by dividing each basis vector by its norm. In \cref{subsec:gram-schmidt}, we will see that finite coproducts in a pre\nobreakdash-Hilbert \(\inv\)\nobreakdash-category can be orthogonalised using a generalisation of the Gram--Schmidt process. Unfortunately, the analogous generalisation of the normalisation step involves certain square roots that do not necessarily exist. Hence, although \(X\) is not necessarily the apex of an \textit{orthonormal} coproduct of \(X_1, \dots, X_n\), it \textit{is} the apex of an \textit{orthogonal} coproduct (a new notion introduced in \cref{sec:orthogonal-biproducts}) of \(X_1, \dots, X_n\). For completeness, the theory of \textit{orthogonal complements} in \(\inv\)\nobreakdash-categories is revised in \cref{sect:orthogonal-complements}. In pre\nobreakdash-Hilbert \(\inv\)\nobreakdash-categories, they are especially well behaved: closed monomorphisms and their orthogonal complements form orthogonal coproducts.

\subsection{Orthogonal biproducts}
\label{sec:orthogonal-biproducts}
The new notion of \textit{orthogonal biproduct} generalises the pre-existing notion of \textit{orthonormal biproduct}, which was recalled in \cref{sect:orthonormal-biproducts}.

\begin{definition}
\label{def:orthogonal-biproduct}
A coproduct \((X, s_1, \dots, s_n)\) in a \(\inv\)\nobreakdash-category is \textit{orthogonal} if the morphisms \(s_1, \dots, s_n\) are pairwise orthogonal. A biproduct \((X, s_1, r_1, \dots, s_n, r_n)\) in a \(\inv\)\nobreakdash-category is \textit{orthogonal} if each of the idempotents \(s_k r_k\) is Hermitian.
\end{definition}

Similar to the relationship between orthonormal coproducts and orthonormal biproducts, orthogonal coproducts and orthogonal biproducts are also really just two different perspectives on a single concept. Seeing this requires some work.

\begin{proposition}
\label{prop:orthogonal-coproduct-biproduct}
A coproduct in a \(\inv\)\nobreakdash-category is orthogonal if and only if it can be extended to an orthogonal biproduct.
\end{proposition}

\begin{proof}
Let \((X, s_1, \dots, s_n)\) be a coproduct of objects \(X_1, \dots, X_n\). First, suppose that this coproduct is orthogonal, that is, that \({s_k}^{\!\inv} \! s_j = 0\) when \(j \neq k\).

Fix \(k \in \set{1, \dots, n}\). By universality of the coproduct, there is a unique morphism \(r_k \colon X \to X_k\) such that \(r_k s_k = 1\) and \(r_k s_j = 0\) for all \(j \neq k\). Observe that 
\[
    {s_\ell}^\inv (s_k r_k)^\inv \! s_j
    = (r_k s_\ell)^\inv {s_k}^{\!\inv} \! s_j
    = \begin{cases}
        {s_k}^{\!\inv} \! s_k
        &\text{if } j = k = \ell,\\
        0&\text{otherwise,}
    \end{cases}
    \quad
    = {s_\ell}^\inv \! s_k r_k s_j,
\]
for all \(j, \ell \in \set{1, \dots, n}\). By universality of the coproduct \((X, s_1, \dots, s_n)\) and the product \((X, {s_1}^{\!\inv}\!, \dots, {s_n}^\inv)\), it follows that the idempotent \(s_k r_k\) is Hermitian, and thus, by \cref{prop:closed-mono-retraction}, that \(s_k\) is a closed monomorphism and \(r_k = ({s_k}^{\!\inv} \! s_k)^{-1} {s_k}^{\!\inv}\).

The tuple \((X, r_1, \dots, r_n)\) is a product because its projections are constructed by composing product projections, namely \({s_1}^{\!\inv}\!, \dots, {s_n}^\inv\), with isomorphisms, namely \(({s_1}^{\!\inv} \! s_1)^{-1}, \dots, ({s_n}^\inv \! s_n)^{-1}\). Combining everything, we see that \((X, s_1, r_1, \dots, s_n, r_n)\) is an orthogonal biproduct.

Conversely, suppose that there are morphisms \(r_k \colon X \to X_k\) for each \(k \in \set{1, \dots, n}\) such that \((X, s_1, r_1, \dots, s_n, r_n)\) is an orthogonal biproduct. By \cref{prop:closed-mono-retraction}, each injection \(s_k\) is again a closed monomorphism and \(r_k = ({s_k}^{\!\inv} \! s_k)^{-1} {s_k}^{\!\inv}\), so
\[{s_k}^{\!\inv} \! s_j = {s_k}^{\!\inv} \! s_k ({s_k}^{\!\inv} \! s_k)^{-1} {s_k}^{\!\inv} \! s_j= {s_k}^{\!\inv} \! s_k r_k s_j = {s_k}^{\!\inv} \! s_k0 = 0\]
for all \(j \neq k\). Hence the coproduct \((X, s_1, \dots, s_n)\) is orthogonal.
\end{proof}

For all \(f \colon X_1 \oplus \dots \oplus X_m \to Y_1 \oplus \dots \oplus Y_n\), if the biproducts are orthogonal then
\begin{equation}
    \label{eq:adjoint-matrix-orthogonal}
    p_j f^\inv i_k
    = ({i_j}^{\!\inv} i_j)^{-1} {i_j}^{\!\inv} f^{\inv} {p_k}^{\inv} (p_k {p_k}^{\!\inv})^{-1}
    = ({i_j}^{\!\inv} i_j)^{-1} (p_k f i_j)^\inv {i_k}^{\!\inv} i_k
\end{equation}
for all \(j \in \set{1, \dots, m}\) and all \(k \in \set{1, \dots, n}\).

\subsection{Orthogonal complements}
\label{sect:orthogonal-complements}
This subsection recounts the well-known theory of orthogonal complements in \(\inv\)\nobreakdash-categories with isometric kernels~\cite{heunen:2010:quantum-logic-dagger-kernel}. More precisely, a slight generalisation of this theory to \(\inv\)\nobreakdash-categories with closed kernels is given.

In a \(\inv\)\nobreakdash-category, an \textit{orthogonal complement} of a monomorphism \(m \colon A \to X\) is a kernel of \(m^\inv\). We write \(m^\perp \colon X \ominus A \to X\) for a chosen orthogonal complement of~\(m\). If \(m\) has an isometric orthogonal complement, then \(m^\perp\) is assumed to be isometric.

The concept of orthogonal complements is of course lifted from the theory of Hilbert spaces. Given a closed subspace \(A\) of a Hilbert space \(X\), the set
\[X \ominus A = \setb[\big]{x \in X}{\innerProd{a}{x} = 0 \text{ for all } a \in A}\]
is also a closed subspace of \(X\), and the inclusion map \(X \ominus A \hookrightarrow X\) is an orthogonal complement in \(\Hilb\) of the inclusion map \(A \hookrightarrow X\). 

Orthogonal complements satisfy the relations
\[
    m \leq m^{\perp\perp},
    \hquad    
    m^{\perp\perp\perp} \cong m^\perp,
    \hquad
    0^\perp \cong 1,
    \hquad
    1^\perp \cong 0,
    \hquad\text{and}\hquad
    m \leq n^\perp \iff n \leq m^\perp,
\]
in the preorder of monomorphisms into a fixed object~(see~\cite[Lemma~1]{heunen:2010:quantum-logic-dagger-kernel}). These relations are reinterpretations of well-known properties of kernels and cokernels.

In many \(\inv\)\nobreakdash-categories, if \(s \colon A \to X\) is a closed monomorphism with an orthogonal complement, then the cospan \((X, s, s^\perp)\) is an orthogonal coproduct of \(A\) and \(X \ominus A\). For example, in \(\inv\)\nobreakdash-categories enriched in abelian groups, this follows from the \textit{splitting lemma}.\footnote{See, e.g.,~\cite[Proposition~1.8.7]{borceux:1994:handbook-categorical-algebra-2} for a proof that works for categories enriched in abelian groups.} Indeed, the morphisms \(s^\perp \colon X \ominus A \to X\) and \(s^\inv \colon X \to A\) form a split short exact sequence: \(s^\perp\) is a kernel of \(s^\inv\), and \(s^\inv\) has a section, namely \(s(s^\inv s)^{-1}\). By the splitting lemma, the cospan \((X, s(s^\inv s)^{-1}, s^\perp)\) is a coproduct. Since \(s\) and \(s(s^\inv s)^{-1}\) differ only by precomposition with an isomorphism, the orthogonal cospan \((X, s, s^\perp)\) is also a coproduct.


Recall, from \cref{s:range}, that each morphism \(f\) in a pre\nobreakdash-Hilbert \(\inv\)\nobreakdash-category has a range factorisation \(f = me\) where \(e\) is coisometric. Since \(e\) is epic, the morphism \(m^\perp\) is a kernel of \(f^\inv\), that is, the familiar equation \((\Ran f)^\perp = \Ker f^\inv\) holds. As noted above, \(\Ran f \leq (\Ran f)^{\perp\perp}\); let \(d\) be the unique morphism such that \(m = m^{\perp\perp}d\). Then \(d\) is both monic and epic~\cite[Corollary~1.1.5]{schneiders:1999:quasi-abelian-sheaves}. In summary, \(f = m^{\perp\perp}de\), where \(m^{\perp\perp}\) is isometric, \(e\) is coisometric, and \(d\) is both monic and epic.

\subsection{The Gram--Schmidt process}
\label{subsec:gram-schmidt}
The Gram--Schmidt process transforms a linearly independent family \((x_1, \dots, x_n)\) of vectors in a Hilbert space \(X\) into an orthogonal family  \((e_1, \dots, e_n)\) of vectors in~\(X\) with the same linear span. The transformed vectors are defined recursively by
\begin{equation}
    e_1 = x_1 \qquad\text{and}\qquad e_{m + 1} = x_{m + 1} - \sum_{k = 1}^m e_k \innerProd{e_k}{e_k}^{-1}\innerProd{e_k}{x_{m + 1}}.
\end{equation}

Instead of sets of vectors in a Hilbert space, we now consider families of morphisms in a pre\nobreakdash-Hilbert \(\inv\)\nobreakdash-category with a common codomain. Such a family is called a \textit{wide cospan}, its elements are called its \textit{legs}, and its common codomain is called its \textit{apex}. For our purposes, the right generalisation of the notion of linearly independent set of vectors is the notion of \textit{split} wide cospan, which is defined as follows.

\begin{definition} A \textit{retraction} of a wide cospan \((s_k \colon A_k \to X)_{k = 1}^n\) is a wide span \((r_k \colon X \to A_k)_{k = 1}^n\) such that, for all \(j, k \in \set{1, 2,\dots, n}\),
\[
    r_k s_j = \begin{cases} 1 &\text{if \(k = j\), and}\\0 &\text{otherwise.}\end{cases}
\]
A wide cospan is \textit{split} if it has a retraction.
\end{definition}

Recall that the \textit{union} of a family of subobjects of an object \(X\) is, if it exists, their supremum in the poset of subobjects of \(X\). The generalisation of linear span for a wide cospan of monomorphisms is, if it exists, the union of the subobjects that its legs represent. In a category that has finite biproducts, a cospan \((s_k \colon A_k \to X)_{k = 1}^n\) is split if and only if the morphism
\(
    \bsmallmat{ s_1 & s_2 & \cdots & s_n } \colon A_1 \oplus A_2 \oplus \dots \oplus A_n \to X
\)
is a split monomorphism, in which case the subobject of \(X\) represented by \(\bsmallmat{ s_1 & s_2 & \cdots & s_n }\) is the union \(s_1 \cup s_2 \cup \dots \cup s_n\) of the subobjects of \(X\) represented by the legs of the cospan.

A wide cospan in a \(\inv\)\nobreakdash-category is \textit{orthogonal} if its legs are pairwise orthogonal.

\begin{proposition}[Gram--Schmidt process]
    \label{prop:gram-schmidt}
    For each split cospan \((s_k \colon A_k \to X)_{k = 1}^n\) in a pre\nobreakdash-Hilbert \(\inv\)\nobreakdash-category, the equations
\begin{equation}
    \label{eq:gram-recursive}
    t_1 = s_1 \qquad\text{and}\qquad t_{m + 1} = s_{m + 1} - \sum_{k = 1}^{m} t_k({t_k}^{\!\inv} t_k)^{-1} {t_k}^{\!\inv} s_{m + 1}
\end{equation}
recursively define an orthogonal cospan \((t_k \colon A_k \to X)_{k = 1}^n\) of closed monomorphisms such that
\(t_1 \cup \dots \cup t_k = s_1 \cup \dots \cup s_k\)
for each \(k \in \set{1, 2, \dots, n}\).
\end{proposition}

Like in the Gram--Schmidt process for Hilbert spaces, the morphism \(t_{m + 1}\) is obtained by orthogonally projecting \(s_{m + 1}\) onto \((t_1 \cup \dots \cup t_m)^\perp\).

\begin{proof}
We proceed by induction on the positive integer \(n\).

For the base case, suppose that \(n = 1\). The cospan \((t_k \colon A_k \to X)_{k = 1}^n\) is vacuously orthogonal. Also, its only leg \(t_1 = s_1\) is split monic, and thus closed monic by \cref{lem:split-is-closed}. Finally, the statement about unions holds because \(t_1 = s_1\).

For the inductive step, suppose that \(n = m + 1\) for some positive integer \(m\), that the cospan \((t_k \colon A_k \to X)_{k = 1}^m\) is orthogonal, that its legs are closed monomorphisms, and that \(t_1 \cup \dots \cup t_k = s_1 \cup \dots \cup s_k\) for all \(k \in \set{1, \dots, m}\).

By the earlier discussion about unions of split wide cospans, there exists an isomorphism
\(u \colon A_1 \oplus \cdots \oplus A_m \to A_1 \oplus \cdots \oplus A_m\)
such that
\begin{equation}
    \label{eq:gram-hypothesis-iso}
    \begin{bmatrix}t_1 & \cdots & t_m \end{bmatrix} = \begin{bmatrix}s_1 & \cdots & s_m \end{bmatrix}u.
\end{equation}
By equations \cref{eq:gram-hypothesis-iso,eq:gram-recursive},
\begin{equation}
    \label{eq:gram-inductive-iso}
    \begin{bmatrix}
        \begin{bmatrix}t_1 & \cdots & t_m
        \end{bmatrix}
        & t_{m + 1}
    \end{bmatrix}
    =
    \begin{bmatrix}
        \begin{bmatrix}s_1 & \cdots & s_m \end{bmatrix}
        & s_{m + 1}
    \end{bmatrix}
    v
\end{equation}
where \(v \colon (A_1 \oplus \cdots \oplus A_m) \oplus A_{m + 1} \to (A_1 \oplus \cdots \oplus A_m) \oplus A_{m + 1}\) is defined by
\[
    v =
    \begin{bmatrix}
        u
        &
        -u
        \begin{bmatrix}
            ({t_1}^{\!\inv} t_1)^{-1}{t_1}^{\!\inv} \\
            \vdots \\
            ({t_m}^\inv t_m)^{-1}{t_m}^\inv
        \end{bmatrix}
        s_{m + 1}
        \\ 0
        & 1
    \end{bmatrix}.
\]

By assumption, the cospan \((s_k \colon A_k \to X)_{k = 1}^{m + 1}\) is split, that is, it has a retraction \((r_k \colon X \to A_k)_{k = 1}^{m + 1}\). Now
\begin{align*}
    r_{m + 1}t_{m + 1}
    &=
    r_{m + 1}
    \begin{bmatrix}
        \begin{bmatrix}t_1 & \cdots & t_m
        \end{bmatrix}
        & t_{m + 1}
    \end{bmatrix}
    i_2
    \\&= r_{m + 1}
    \begin{bmatrix}
        \begin{bmatrix}s_1 & \cdots & s_m \end{bmatrix}
        & s_{m + 1}
    \end{bmatrix}
    v i_2
    = p_2 v i_2
    = 1,
\end{align*}
so \(r_{m + 1}\) is a retraction of \(t_{m + 1}\). By \cref{lem:split-is-closed}, \(t_{m + 1}\) is a closed monomorphism.

For each \(j \in \set{1, 2, \dots, m}\), the morphism \(t_j\) is orthogonal to \(t_{m + 1}\); indeed,
\begin{align*}
{t_j}^{\!\inv} t_{m + 1}
    &= {t_j}^{\!\inv} \! s_{m + 1} - {t_j}^{\!\inv} \sum_{k = 1}^m t_k({t_k}^{\!\inv} t_k)^{-1} {t_k}^{\!\inv}  s_{m + 1}
    \\&\qquad\qquad= {t_j}^{\!\inv} \! s_{m + 1} - {t_j}^{\!\inv} t_j({t_j}^{\!\inv} t_j)^{-1} {t_j}^{\!\inv}  s_{m + 1}
    = 0
\end{align*}
because \(t_j\) is orthogonal to each \(t_k\) with \(k \in \set{1, \dots, m}\) and \(k \neq j\).

So far we have shown that the cospan \((t_k \colon A_k \to X)_{k = 1}^{m + 1}\) is orthogonal and that its legs are closed monomorphisms. Hence, for all \(j, k \in \set{1, 2, \dots, m + 1}\),
\[
    ({t_k}^{\!\inv} t_k)^{-1} {t_k}^{\!\inv} t_j = \begin{cases}
        1 &\text{if \(k = j\), and}\\
        0 &\text{otherwise.}
    \end{cases}
\]
This means that the cospan is split, so the union \(t_1 \cup \dots \cup t_m \cup t_{m + 1}\) exists and is represented by the morphism \(\bsmallmat{
    \bsmallmat{t_1 & \cdots & t_m
    }
    & t_{m + 1}
}\). The union \(s_1 \cup \dots \cup s_m \cup s_{m + 1}\) is similarly represented by \(\bsmallmat{
    \bsmallmat{s_1 & \cdots & s_m }
    & s_{m + 1}
}\). Observe that \(v\) is an isomorphism with
\[v^{-1} = \begin{bmatrix}
    u^{-1}
    &
    \begin{bmatrix}({t_1}^{\!\inv} t_1)^{-1}{t_1}^{\!\inv} \\ \vdots \\ ({t_m}^\inv t_m)^{-1}{t_m}^\inv
    \end{bmatrix}s_{m + 1}
    \\ 0
    & 1
\end{bmatrix}.\]
The equality \(t_1 \cup \dots \cup t_{m + 1} = s_1 \cup \dots \cup s_{m + 1}\) thus follows from equation \cref{eq:gram-inductive-iso}.
\end{proof}

The Gram--Schmidt process transforms coproducts into orthogonal coproducts. Indeed, consider a coproduct \((X, s_1, \dots, s_n)\) of objects \(A_1, \dots, A_n\) in a pre\nobreakdash-Hilbert \(\inv\)\nobreakdash-category. The morphism
\(\bsmallmat{s_1 & \cdots & s_n}\)
is invertible, so the cospan \((s_k \colon A_k \to X)_{k = 1}^n\) is split. The Gram--Schmidt process yields an orthogonal cospan \((t_k \colon A_k \to X)_{k = 1}^n\) of closed monomorphisms such that \(t_1 \cup \dots \cup t_n = s_1 \cup \dots \cup s_n\). As \(\bsmallmat{s_1 & s_2 & \cdots & s_n}\) is invertible, the morphism \(\bsmallmat{t_1 & \cdots & t_n}\) is also invertible, so \((X, t_1, \dots, t_n)\) is an orthogonal coproduct of \(A_1, \dots, A_n\). If the morphisms \(({t_1}^{\!\inv} t_1)^{-1}, \dots, ({t_n}^{\!\inv} t_n)^{-1}\) respectively have Hermitian square roots \(r_1, \dots r_n\), then \((X,\, t_1 r_1, \, \dots,\, t_nr_n)\) is an orthonormal coproduct of \(A_1, \dots, A_n\).

\section{Order and inner products}
\label{sect:examples}

To justify introducing an abstraction like pre\nobreakdash-Hilbert \(\inv\)\nobreakdash-category, it ought to have many interesting examples. This section describes two families of examples of pre\nobreakdash-Hilbert \(\inv\)\nobreakdash-categories. The first, which is the topic of \cref{sect:fin-inner-prod-modules}, consists of the \(\inv\)\nobreakdash-categories \(\FinInnerProd{D}\) of finite-dimensional inner product \(D\)\nobreakdash-modules and \(D\)-linear maps, where \(D\) is an ordered division \(\inv\)\nobreakdash-ring. The second---the topic of \cref{sect:hilbert-W-modules}---consists of the \(\inv\)\nobreakdash-categories \(\PHilb{A}\) of self-dual Hilbert \(A\)\nobreakdash-modules and bounded \(A\)-linear maps, where \(A\) is a complex W*-algebra. Setting \(A = \Comps\) recovers the \(\inv\)\nobreakdash-category \(\Hilb\) of complex Hilbert spaces that was introduced in \cref{sect:background}. The concepts of \textit{ordered \(\inv\)\nobreakdash-ring} and \textit{inner product module}, which are foundational for both families of examples, are introduced in \cref{sect:preliminaries}.
In fact all pre\nobreakdash-Hilbert \(\inv\)\nobreakdash-categories are, in some sense, \(\inv\)\nobreakdash-categories of inner product modules; this claim is made precise in \cref{prop:canonical-structure} of \cref{sect:r-rings}.

\subsection{Ordered \texorpdfstring{\(\inv\)}{*}-rings and their inner product modules}
\label{sect:preliminaries}
A \textit{\(\inv\)\nobreakdash-ring} is a ring \(R\) together with a function \(\blank^\inv \colon R \to R\) such that, for all \(r, s \in R\),
\[
    1^\inv = 1, \qquad
    (sr)^\inv = r^\inv \! s^\inv\!, \qquad
    (r^\inv)^\inv = r, \qquad\text{and}\qquad
    (r + s)^\inv = r^\inv + s^\inv.
\]
The operation \(r \mapsto r^\inv\) is called the \textit{involution} of the \(\inv\)\nobreakdash-ring. A \(\inv\)\nobreakdash-ring will be called \textit{anisotropic} if \(r^\inv r = 0\) implies \(r = 0\) for each \(r \in R\).\footnote{The term \textit{anisotropic} comes from the literature on Hermitian forms~(see, e.g., \cite{knus:1991:quadratic-hermitian-forms-rings}). In the literature on \(\inv\)\nobreakdash-rings (see, e.g.,~\cite{berberian:baer-star-rings}), the term \textit{proper} is used instead.}

An element \(a\) of a \(\inv\)\nobreakdash-ring is called \textit{Hermitian} if \(a^\inv = a\). For each \(\inv\)\nobreakdash-ring \(R\), the set of Hermitian elements of \(R\) will be denoted \(\SelfAdj{R}\). If \(a, b \in \SelfAdj{R}\) and \(r \in R\), then \(a^\inv \in \SelfAdj{R}\) and \(r^\inv a r + b \in \SelfAdj{R}\). The set \(\SelfAdj{R}\) is not, in general, closed under multiplication.

\begin{definition}
An \textit{ordering} of a \(\inv\)\nobreakdash-ring \(R\) is a partial order \(\leq\) on \(\SelfAdj{R}\) such that
\[
    0 \leq 1,
    \qquad\qquad\text{and}\qquad\qquad
    a \leq b \;\implies\; r^\inv \!ar + c \leq r^\inv br + c,
\]
for all \(a, b, c \in \SelfAdj{R}\) and all \(r \in R\).
\end{definition}

We sometimes write \(a \geq b\) instead of \(b \leq a\).

\begin{definition}
A \textit{positive cone} on a \(\inv\)\nobreakdash-ring \(R\) is a subset \(P\) of \(\SelfAdj{R}\) such that
\[
    1 \in P,
    \qquad
    P + P \subseteq P,
    \qquad
    r^\inv P r \subseteq P,
    \qquad\text{and}\qquad
    P \cap -P = \set{0}
\]
for all \(r \in R\).
\end{definition}

\begin{proposition}
For each \(\inv\)\nobreakdash-ring \(R\), the equations
\[P = \setb[\big]{a \in \SelfAdj{R}}{0 \leq a} \qquad\text{and}\qquad a \leq b \;\iff\; b - a \in P\]
define a bijection between the orderings \(\leq\) of \(R\) and the positive cones \(P\) on \(R\).
\end{proposition}

\begin{proof}
First, let \(\leq\) be an ordering of \(R\), and define \(P\) in terms of \(\leq\) as above. Then \(1 \in P\) because \(0 \leq 1\). Also, if \(a, b \in P\) and \(r \in R\), then
\[0 \leq b = r^\inv 0 r + b \leq r^\inv a r + b,\]
and so \(r^\inv a r + b \in P\) by transitivity of \(\leq\). Finally, if \(a \in P \cap -P\), then \(a \in P\) and \(a = -b\) for some \(b \in P\), and so
\[0 \leq a = a + 0 \leq a + b = 0,\]
which means that \(a = 0\) by antisymmetry of \(\leq\). Hence \(P\) is a positive cone on \(R\).

Conversely, let \(P\) be a positive cone on \(R\), and define \(\leq\) in terms of \(P\) as above. For reflexivity of \(\leq\), for all \(a \in \SelfAdj{R}\), we have \(a \leq a\) because
\[a - a = 0 \in P \cap -P \subseteq P.\]
For antisymmetry of \(\leq\), for all \(a, b \in \SelfAdj{R}\), if \(a \leq b\) and \(b \leq a\), then \(a - b \in P\) and \(b - a \in P\), so \(a - b \in P \cap -P = \set{0}\), and thus \(a - b = 0\), that is \(a = b\). For transitivity of \(\leq\), for all \(a, b, c \in \SelfAdj{R}\), if \(a \leq b\) and \(b \leq c\) then \(b - a \in P\) and \(c - b \in P\), so \[c - a = (c - b) + (b - a) \in P + P \subseteq P.\]
Also \(0 \leq 1\) because \(1 \in P\), and, for all \(a, b, c \in \SelfAdj{R}\) and all \(r \in R\), if \(a \leq b\) then \(r^\inv \! a r + c \leq r^\inv b r + c\) because
\[(r^\inv b r + c) - (r^\inv \! a r + c) = r^\inv(b - a)r \in r^\inv P r \subseteq P. \qedhere\]
\end{proof}

\begin{definition}
An \textit{ordered \(\inv\)\nobreakdash-ring} is an anisotropic \(\inv\)\nobreakdash-ring \(R\) equipped with an ordering \(\leq\), or equivalently, with a positive cone \(\Pos{R}\). The elements of \(\Pos{R}\) are called the \textit{positive} elements of \(R\). A \textit{totally ordered \(\inv\)\nobreakdash-ring} is an ordered \(\inv\)\nobreakdash-ring in which the ordering \(\leq\) is total, or equivalently, in which \(\Pos{R} \cup -\Pos{R} = \SelfAdj{R}\).
\end{definition}

Ordered \(\inv\)\nobreakdash-rings are the \(\inv\)\nobreakdash-ring analogue of \citeauthor{schotz:equivalence-order-algebraic}'s \textit{ordered \(\inv\)\nobreakdash-algebras}~\cite{schotz:equivalence-order-algebraic}. Also, positive cones have almost the same axioms as \citeauthor{cimpric:2009:quadratic-module}'s \textit{quadratic modules}~\cite{cimpric:2009:quadratic-module}. In fact, every positive cone on a non-zero \(\inv\)\nobreakdash-ring is a quadratic module. Total orderings are a generalisation of \citeauthor{prestel:1984:lectures-on-formally-real}'s \textit{semiorderings} of fields~\cite{prestel:1984:lectures-on-formally-real} and of \citeauthor{baer:1952:chapter-dualities}'s \textit{domains of positivity} of division \(\inv\)\nobreakdash-rings~\cite[Appendix~I]{baer:1952:chapter-dualities}, which are now called \textit{Baer orderings} (see, e.g.,~\cite{holland:1977:orderings-and-square-roots,craven:1995:orderings-valuations-hermitian}).

\begin{definition}
    Let \(R\) be an ordered \(\inv\)\nobreakdash-ring, and let \(X\) be a right \(R\)\nobreakdash-module. An \textit{inner product} on~\(X\) is a function \(\innerProd{\blank}{\blank} \colon X \times X \to R\) with the following properties:
    \begin{itemize}
        \item \textit{linearity}: \(\innerProd{x}{yr + z} = \innerProd{x}{y}r + \innerProd{x}{z}\) for all \(x, y, z \in X\) and all \(r \in R\);
        \item \textit{Hermitianness}: \(\innerProd{x}{y}^\inv = \innerProd{y}{x}\) for all \(x, y \in X\);
        \item \textit{positivity}: \(\innerProd{x}{x} \geq 0\) for all \(x \in X\); and
        \item \textit{anisotropy}: \(\innerProd{x}{x} = 0\) implies \(x = 0\) for all \(x \in X\).
    \end{itemize}
    An \textit{inner product \(R\)\nobreakdash-module} is a right \(R\)\nobreakdash-module equipped with an inner product.
\end{definition}

Inner product modules over totally ordered division \(\inv\)\nobreakdash-rings were considered by~\textcite{holland:star-valuations}.

Let \(X\) and \(Y\) be inner product \(R\)\nobreakdash-modules. For each function \(f \colon X \to Y\), there is at most one function \(f^\inv \colon Y \to X\) such that \(\innerProd{f^\inv y}{x} = \innerProd{y}{fx}\) for all \(x \in X\) and all \(y \in Y\). If such a function \(f^\inv\) exists, then \(f\) is said to be \textit{adjointable} and \(f^\inv\) is called the \textit{adjoint} of \(f\). Every adjointable function \(f \colon X \to Y\) is \(R\)-linear.

Inner product \(R\)\nobreakdash-modules and adjointable functions form a \(\inv\)\nobreakdash-category, which we denote by \(\InnerProd{R}\). It satisfies almost all of the pre\nobreakdash-Hilbert \(\inv\)\nobreakdash-category axioms, inheriting an enrichment in abelian groups and finite orthonormal biproducts from the category \(\Mod{R}\) of right \(R\)\nobreakdash-modules and \(R\)-linear maps.

\begin{proposition}
\label{prop:inner-prod-axioms}
The \(\inv\)\nobreakdash-category \(\InnerProd{R}\) satisfies axioms \cref{axiom:R-zero,axiom:R-coproduct,axiom:R-normal}.
\end{proposition}

\begin{proof}
The zero module \(\zero = \set{0}\) is a zero object in \(\InnerProd{R}\) when equipped with its unique inner product. Hence axiom \cref{axiom:R-zero} holds.

Recall that addition of morphisms in \(\Mod{R}\) is defined pointwise in terms of addition in their codomain. For all morphisms \(f, g \colon X \to Y\) of \(\InnerProd{R}\), the morphisms \(f + g\) and \(-f\) of \(\Mod{R}\) are also morphisms of \(\InnerProd{R}\). Indeed,
\begin{align*}
    \innerProd{y}{(f + g)x}
    &= \innerProd{y}{fx + gx}
    = \innerProd{y}{fx} + \innerProd{y}{gx}
    \\&= \innerProd{f^\inv y}{x} + \innerProd{g^\inv y}{x}
    = \innerProd{f^\inv y + g^\inv y}{x}
    = \innerProd{(f^\inv + g^\inv)y}{x}
\end{align*}
and
\[
    \innerProd{y}{(-f)x} = \innerProd{y}{-fx} = -\innerProd{y}{fx} = -\innerProd{f^\inv y}{x} = \innerProd{-f^\inv y}{x} = \innerProd{(-f^\inv)y}{x}
\]
for all \(x \in X\) and all \(y \in Y\), so \((f + g)^\inv = f^\inv + g^\inv\) and \((-f)^\inv = -f^\inv\). It follows that \(\InnerProd{R}\) inherits the enrichment of \(\Mod{R}\) in abelian groups.

Also recall that each pair of objects of \(\Mod{R}\) has a canonical biproduct: their direct sum together with the canonical projections and injections. For all objects \(X_1\) and \(X_2\) of \(\InnerProd{R}\), and all biproducts \((X, s_1, r_1, s_2, r_2)\) of \(X_1\) and \(X_2\) in~\(\Mod{R}\), there is a unique inner product on \(X\) such that \({s_1}^{\!\inv} = r_1\) and \({s_2}^{\!\inv} = r_2\). Indeed, if such an inner product exists, then it is unique because, for all \(x, y \in X\),
\[
    \innerProd{x}{y}
    = \innerProd[\big]{x}{(s_1r_1 + s_2r_2)y}
    = \innerProd{x}{s_1r_1y} + \innerProd{x}{s_2r_2y}
    = \innerProd{r_1x}{r_1y} + \innerProd{r_2x}{r_2y}.
\]
For existence, define the inner product using this equation. Hermitianness, positivity and linearity are easy to check. For anisotropy, if \(\innerProd{x}{x} = 0\), then
\[0 \leq \innerProd{r_1x}{r_1x} \leq \innerProd{r_1x}{r_1x} + \innerProd{r_2x}{r_2x} = \innerProd{x}{x} = 0,\]
so \(\innerProd{r_1x}{r_1x} = 0\) by antisymmetry, and thus \(r_1x = 0\) by anisotropy; similarly \(r_2x = 0\), so \(x = 0\). Finally \({s_1}^{\!\inv} = r_1\) because, for all \(x \in X\) and all \(y_1 \in X_1\),
\[\innerProd{x}{s_1 y_1} = \innerProd{r_1x}{r_1s_1y_1} + \innerProd{r_2x}{r_2s_1y_1} = \innerProd{r_1x}{y_1}  + \innerProd{r_2x}{0} = \innerProd{r_1x}{y_1},\]
and \({s_2}^{\!\inv} = r_2\) similarly. With this inner product, the tuple \((X, s_1, r_1, s_2, r_2)\) is an orthonormal biproduct of \(X_1\) and \(X_2\) in \(\InnerProd{R}\). Hence axiom \cref{axiom:R-coproduct} holds.

Finally, \(\diagonal \colon X \to X \oplus X\) is a kernel of \(\copair{1}{-1}\), so axiom \cref{axiom:R-normal} holds.
\end{proof}

Unfortunately, for some \(R\), there are morphisms in \(\InnerProd{R}\) whose kernel in \(\Mod{R}\) is not adjointable, and, even worse, that do not have an isometric kernel at all (see \cref{f:no-isometric-kernel}). In other words, the \(\inv\)\nobreakdash-category \(\InnerProd{R}\) does not, in general, satisfy axiom~\cref{axiom:R-kernel}. However, if~\(R\) is nice enough, certain full subcategories of it are pre\nobreakdash-Hilbert \(\inv\)\nobreakdash-categories. The examples in \cref{sect:fin-inner-prod-modules,sect:hilbert-W-modules} are of this form.

\subsection{Finite-dimensional inner product modules}
\label{sect:fin-inner-prod-modules}
A \textit{division ring} is a non-zero ring in which every non-zero element is invertible. Every module over a division ring is free, and all bases of such a module have the same cardinality, called its \textit{dimension}.

The rational numbers \(\Rats\), the real numbers \(\Reals\), the complex numbers \(\Comps\), and the quaternions \(\Quats\) are all canonically totally ordered division \(\inv\)\nobreakdash-rings. A more exotic example~\cite[Example~2.10]{craven:1995:orderings-valuations-hermitian} is the field \(\Reals(\!(X)\!)\) of formal Laurent series over \(\Reals\), where
\[\paren[\Big]{\sum_{k = n}^\infty a_k X^k}^\inv = \sum_{k = n}^\infty a_k (-1)^k X^k,\]
and
\[\Pos{\Reals(\!(X)\!)} = \setb[\Big]{\sum_{k = n}^\infty a_k X^{2k} \in \Reals(\!(X)\!)}{(-1)^n a_{n} > 0} \cup \set[\big]{0}.\]

Let \(D\) be an ordered division \(\inv\)\nobreakdash-ring, and let \(\FinInnerProd{D}\) be the category of finite-dimensional inner product \(D\)\nobreakdash-modules and \(D\)-linear maps.

\begin{proposition}
\label{p:finn-prod-r-cat}
The category \(\FinInnerProd{D}\) is a full subcategory of \(\InnerProd{D}\). With the inherited \(\inv\)\nobreakdash-category structure, it is actually a pre\nobreakdash-Hilbert \(\inv\)\nobreakdash-category.
\end{proposition}

\begin{proof}
The well-known theory of finite-dimensional inner product spaces over \(\Reals\) and~\(\Comps\) generalises, \textit{mutatis mutandis}, to finite-dimensional inner product \(D\)\nobreakdash-modules. In particular, this is true of the usual construction of adjoints:
\begin{itemize}
    \item First, bases can be orthogonalised using the \textit{Gram--Schmidt process}. Indeed, if \(X\) is an inner product \(D\)\nobreakdash-module with basis \(\set{x_1, x_2, \dots, x_n}\), then the equations
    \[e_1 = x_1 \qquad\text{and}\qquad e_{m + 1} = x_{m + 1} - \sum_{k = 1}^m e_k \innerProd{e_k}{e_k}^{-1}\innerProd{e_k}{x_{m+1}}\]
    recursively define an orthogonal basis \(\set{e_1, e_2, \dots, e_n}\) of \(X\).
    \item Next, \textit{Parseval's identity}
    \[x = \sum_{k = 1}^n e_k \innerProd{e_k}{e_k}^{-1} \innerProd{e_k}{x}\]
    holds for each \(x \in X\).
    \item Finally, the map \(g \colon Y \to X\) defined by
    \[gy = \sum_{k = 1}^n e_k \innerProd{e_k}{e_k}^{-1}\innerProd{fe_k}{y}\]
    is the adjoint of \(f \colon X \to Y\). To prove this, use Parseval's identity for \(x\).
\end{itemize}
As all morphisms in \(\FinInnerProd{D}\) are adjointable and all morphisms in \(\InnerProd{D}\) are linear, the category \(\FinInnerProd{D}\) is a full subcategory of \(\InnerProd{D}\), and so it inherits from \(\InnerProd{D}\) the structure of a \(\inv\)\nobreakdash-category and an enrichment in abelian groups. It similarly inherits finite orthonormal biproducts because \(\zero\) is finite dimensional and the direct sum of two finite-dimensional \(D\)\nobreakdash-modules is again finite dimensional. Hence \(\FinInnerProd{D}\) satisfies axioms \cref{axiom:R-zero,axiom:R-coproduct,axiom:R-normal}.

Unlike \(\InnerProd{D}\), the \(\inv\)\nobreakdash-category \(\FinInnerProd{D}\) also satisfies axiom \cref{axiom:R-kernel}. Indeed, its kernels are inherited from \(\Mod{D}\). First, recall that every morphism \(f \colon X \to Y\) of \(\Mod{D}\) has a canonical kernel: the inclusion into \(X\) of the submodule
\[\Ker f = \setb[\big]{x \in X}{fx = 0}\]
of \(X\). 
For each morphism \(f \colon X \to Y\) in \(\FinInnerProd{D}\), and each kernel \(s \colon A \to X\) of~\(f\) in \(\Mod{D}\), we will show that the \(D\)\nobreakdash-module \(A\) is finite dimensional and has a unique inner product that makes \(s\) an isometry in \(\FinInnerProd{D}\). 

First, a finite basis \(\mathcal{B}\) of \(A\) can be constructed as follows: starting with \(\mathcal{B} = \emptyset\), while there is a vector \(a\) in \(A\) that is not in the span of \(\mathcal{B}\), add \(a\) to \(\mathcal{B}\). Throughout the procedure, the subset \(\mathcal{B}\) of \(A\) remains linearly independent. As \(s\) is injective, this is also true of the subset \(s(\mathcal{B})\) of \(X\). Also \(s(\mathcal{B})\) and \(\mathcal{B}\) always have the same cardinality. Since this cardinality increases throughout the procedure, and is bounded above by the dimension of \(X\), the procedure must terminate. When it does, the set \(\mathcal{B}\) spans~\(A\), and so is a basis of \(A\).

Next, if such an inner product on \(A\) exists, then it is unique because
\[\innerProd{a}{b} = \innerProd{s^\inv \! s a}{s^\inv \! s b} = \innerProd{sa}{ss^\inv \! s b} = \innerProd{sa}{sb}\]
for all \(a,b \in A\). For existence, define the inner product on \(A\) by this equation. Then \(s\) is a morphism of \(\FinInnerProd{D}\), and so, by the discussion above, is adjointable. Now \(\innerProd{s^\inv \! s a}{b} = \innerProd{sa}{sb} = \innerProd{a}{b}\)
for all \(a, b \in A\), so \(s^\inv \! s = 1^\inv = 1\).

Since the forgetful functor from \(\FinInnerProd{D}\) to \(\Mod{D}\) is full and faithful, it reflects limits. Hence \(s\) is an isometric kernel of \(f\) in \(\FinInnerProd{D}\).
\end{proof}

It is well known that \(\FinInnerProd{\Comps}\) is equivalent to the \(\inv\)\nobreakdash-category \(\Mat{\Comps}\) of natural numbers and complex matrices. More generally, if every positive element of \(D\) has a positive square root, then \(\FinInnerProd{D}\) is equivalent to the \(\inv\)\nobreakdash-category \(\Mat{D}\) of natural numbers and \(D\)-valued matrices; the square roots are needed to normalise the orthogonal bases produced by the Gram--Schmidt procedure. Without this assumption about square roots, the \(\inv\)\nobreakdash-category \(\FinInnerProd{D}\) is still equivalent to a \(\inv\)\nobreakdash-category whose morphisms are matrices, namely, the \(\inv\)\nobreakdash-category \(\WeightMat{D}\) of \(D\)\nobreakdash-valued weights and matrices, which is defined as follows.

The objects of \(\WeightMat{D}\) are pairs \((m, \alpha)\) where \(m\) is a natural number and \(\alpha\) is a vector of length \(m\) whose entries, called \textit{weights}, are non-zero positive elements of~\(D\). The morphisms \((m, \alpha) \to (n, \beta)\) in \(\WeightMat{D}\) are matrices with \(n\) rows and \(m\) columns whose entries are elements of \(D\); composition is matrix multiplication. The involution of \(\WeightMat{D}\) is defined on each morphism \(M \colon (m, \alpha) \to (n, \beta)\) by
\begin{equation}
\label{eq:weight-mat-inv}
    (M^\inv)_{jk} = \alpha_j {M_{kj}}^{\!\inv} {\beta_k}^{\!-1}
\end{equation}
for all \(j \in \set{1, \dots, m}\) and \(k \in \set{1, \dots, n}\).

The equivalence \(\WeightMat{D} \to \FinInnerProd{D}\) sends an object \((m, \alpha)\) to the \(D\)\nobreakdash-module \(D^m\) equipped with the inner product
\[\innerProd*{x}{y} = {x_1}^{\!\inv} {\alpha_1}^{\!-1} y_1 + \dots + {x_m}^{\!\inv} {\alpha_m}^{\!-1} y_m,\]
and it sends a morphism \(M \colon (m, \alpha) \to (n, \beta)\) to the \(D\)-linear map \(D^m \to D^n\) defined by \(x \mapsto Mx\). It is easy to check that it is a full and faithful \(\inv\)\nobreakdash-functor. For unitary essential surjectivity, each \(n\)-dimensional inner product \(D\)\nobreakdash-module has an orthogonal basis \(\set{e_1, \dots, e_n}\) by the Gram--Schmidt process, and so is unitarily isomorphic to the image of the weight vector \(\paren[\big]{\innerProd{e_1}{e_1}^{-1}, \dots, \innerProd{e_n}{e_n}^{-1}}\) under the equivalence.

Observe that equation \cref{eq:weight-mat-inv} is equivalent to
\[
    (M^\inv)_{jk}\beta_k = (M_{kj} \alpha_j)^\inv,
\]
which is a generalisation of Bayes' law. This equation is a special case of equation~\cref{eq:adjoint-matrix-orthogonal} because weight vectors of length~\(n\) are mapped by the equivalence to orthogonal biproducts of \(n\) copies of \(D\).

\subsection{Hilbert W*-modules}
\label{sect:hilbert-W-modules}
The original articles~\cite{kaplansky:modules-operator-algebras,paschke:inner-product-modules-b-algebras} on Hilbert modules over C*-algebras remain good introductions to the topic. The reader may also find helpful the PhD theses of the Westerbaan brothers~\cite{westerbaan:dagger-and-dilation,westerbaan:category-von-neumann-algebras}, which give a comprehensive account of W*-algebras and their Hilbert modules.

A (\textit{complex unital}) \textit{C*-algebra} is a \(\inv\)\nobreakdash-ring \(A\) equipped with a \(\inv\)\nobreakdash-ring homomorphism from \(\Comps\) to the centre of \(A\), and a complete norm \(\norm{\blank} \colon A \to \PosReals\), such that
\[\norm{ab} \leq \norm{a}\norm{b} \qquad\text{and}\qquad \norm{a^\inv a} = \norm{a}^2\]
for all \(a, b \in A\). Every C*-algebra \(A\) is an ordered \(\inv\)\nobreakdash-ring when equipped with the positive cone \(\Pos{A} = \setb{a^\inv a}{a \in A}\).

Let \(A\) be a C*-algebra. Every inner product \(A\)\nobreakdash-module is canonically normed, with \(\norm{x} = \sqrt{\norm{\innerProd{x}{x}}}\). A \textit{Hilbert \(A\)\nobreakdash-module} (or \textit{Hilbert C*-module over \(A\)}) is an inner product \(A\)\nobreakdash-module that is complete with respect to its canonical norm. A map \(f \colon X \to Y\) between inner product \(A\)\nobreakdash-modules is \textit{bounded} if there is a \(c \in \Pos{\Reals}\) such that \(\norm{fx} \leq c \norm{x}\) for all \(x \in X\). An inner product \(A\)\nobreakdash-module \(X\) is \textit{self-dual} if, for each bounded \(A\)-linear map \(f \colon X \to A\), there exists \(u \in X\) such that \(fx = \innerProd{u}{x}\) for all \(x \in X\). All self-dual inner product \(A\)\nobreakdash-modules are Hilbert \(A\)\nobreakdash-modules~\cite[449]{paschke:inner-product-modules-b-algebras}. Let \(\Hilb_A\) denote the category of self-dual Hilbert \(A\)\nobreakdash-modules and bounded \(A\)-linear maps. This definition is consistent with the earlier definition of \(\Hilb_\Comps\) in \cref{sect:background}.

\begin{proposition}
\label{p:hilb_module_axioms}
The category \(\Hilb_A\) is a full subcategory of \(\InnerProd{A}\). With the inherited \(\inv\)\nobreakdash-category structure, it satisfies axioms \cref{axiom:R-zero,axiom:R-coproduct,axiom:R-normal}.
\end{proposition}

\begin{proof}
A function between self-dual Hilbert \(A\)\nobreakdash-modules is bounded and \(A\)-linear if~\cite[\S~8.1.7]{blecher-merdy:operator-algebras} and only if~\cite[\S~3.4]{paschke:inner-product-modules-b-algebras} it is adjointable. Hence the category \(\Hilb_A\) is a full subcategory of \(\InnerProd{A}\). It thus inherits from \(\InnerProd{A}\) an enrichment in abelian groups and the structure of a \(\inv\)\nobreakdash-category. It similarly inherits finite orthonormal biproducts because the inner product \(A\)\nobreakdash-module \(\zero\) is self-dual and the orthonormal biproduct of two self-dual inner product \(A\)\nobreakdash-modules is again self-dual. For the first claim, the zero map is the only bounded \(A\)-linear map from \(\zero\) to \(A\), and \(0(x) = \innerProd{0}{x}\) for all \(x \in \zero\). For the second claim, let \(X\) and \(Y\) be self-dual Hilbert \(A\)\nobreakdash-modules and let \(f \colon X \oplus Y \to A\) be a bounded \(A\)-linear map. As \(fi_1\) is bounded and \(A\)-linear, and \(X\) is self-dual, there exists \(u \in X\) such that \(fi_1x = \innerProd{u}{x}\). Similarly, there exists \(v \in Y\) such that \(fi_2 y = \innerProd{v}{y}\). It follows that
\[f(x, y) = fi_1x + fi_2y = \innerProd{u}{x} + \innerProd{v}{y} = \innerProd{(u, v)}{(x, y)}. \qedhere\]
\end{proof}

The \(\inv\)\nobreakdash-category \(\Hilb_A\) does not, in general, satisfy axiom \cref{axiom:R-kernel}.

\begin{remark}
\label{f:no-isometric-kernel}
Let \(X\) be a compact Hausdorff space, and let \(C(X)\) be the C*\nobreakdash-algebra of continuous functions \(X \to \Comps\). If \(X\) is not totally disconnected (e.g., if \(X\) is the closed interval \([0, 1] \subseteq \Reals\)), then there is a morphism~\(f\) in \(\Hilb_{C(X)}\) that does not have an isometric kernel (adapt~\cite[Proposition~10.3]{heunen:frobenius-hilbert-modules}). Furthermore, as \(\Hilb_{C(X)}\) is a full subcategory of \(\InnerProd{C(X)}\), and an object of \(\InnerProd{C(X)}\) comes from \(\Hilb_{C(X)}\) if it admits an isometry to an object from \(\Hilb_{C(X)}\), the morphism \(f\) also does not have an isometric kernel when viewed as a morphism in \(\InnerProd{C(X)}\).
\end{remark}

However, the story is different when \(A\) is a \textit{W*-algebra}.\footnote{W*-algebras are also called \textit{von Neumann algebras}.} Self-dual Hilbert modules over a W*-algebra \(A\) are also called \textit{Hilbert W*-modules} over \(A\).

\begin{proposition}
\label{p:w-module-r-cat}
If \(A\) is a W*-algebra, then \(\Hilb_A\) is a pre\nobreakdash-Hilbert \(\inv\)\nobreakdash-category.
\end{proposition}

The only fact about W*-algebras needed for the proof is that a submodule \(M\) of a Hilbert W*\nobreakdash-module \(X\) is self-dual with respect to its inherited inner product if and only if \(M = M^{\perp\perp}\) (see~\cite[Corollary~1]{frank-troitsky:lefschetz-numbers} and~\cite[149V and 160IV]{westerbaan:dagger-and-dilation}), where for each subset \(S\) of \(X\),
\[S^\perp = \setb[\big]{x \in X}{\innerProd{s}{x} = 0 \text{ for all } s \in S}\]
is the \textit{orthogonal complement} of \(S\) in \(X\). Of course, once we know that \(\Hilb_A\) is a pre\nobreakdash-Hilbert \(\inv\)\nobreakdash-category, this notion of orthogonal complement coincides with the abstract one discussed in \cref{sect:orthogonal-complements}.

\begin{proof}[Proof of \cref{p:w-module-r-cat}]
Axioms \cref{axiom:R-zero,axiom:R-coproduct,axiom:R-normal} hold by \cref{p:hilb_module_axioms}. For axiom \cref{axiom:R-kernel}, let \(f \colon X \to Y\) be a morphism in \(\Hilb_A\). Equip the \(A\)-submodule
\[\Ker f = \setb[\big]{x \in X}{fx = 0}\]
of \(X\) with the restriction of the inner product of \(X\).

We wish to show that \(\Ker f\) is self-dual. The following facts are proved in the same way as the analogous well-known results in the theory of Hilbert spaces:
\begin{itemize}
    \item \(\Ker f = (\Ran f^\inv)^\perp\) where
    \(\Ran f^\inv = \setb[\big]{f^\inv y}{y \in Y}\) is the \textit{range} of \(f^\inv\), and
    \item if \(S \subseteq X\) then \(S^{\perp\perp\perp} = S^\perp\).
\end{itemize}
Combining these facts,
\[(\Ker f)^{\perp \perp} = (\Ran f^\inv)^{\perp\perp\perp} = (\Ran f^\inv)^\perp = \Ker f,\]
so \(\Ker f\) is self-dual by the discussion above.

The inclusion \(s \colon \Ker f \hookrightarrow X\) is a bounded \(A\)-linear map between self-dual Hilbert \(A\)\nobreakdash-modules, so it is adjointable. Also \[\innerProd{s^\inv s x}{x'} = \innerProd{sx}{sx'} = \innerProd{x}{x'}\]
for all \(x, x' \in \Ker f\), so \(s^\inv s = 1\), that is, \(s\) is isometric. It is actually an isometric kernel of \(f\). Indeed, let \(g \colon Z \to X\) be a morphism in \(\Hilb_A\) such that \(fg = 0\). Let \(h \colon Z \to \Ker f\) be the unique morphism in \(\Mod{A}\) such that \(sh = g\). Then \(h\) is adjointable with adjoint \(g^\inv s\) because
\[\innerProd{x}{hz} = \innerProd{sx}{shz} = \innerProd{sx}{gz} = \innerProd{g^\inv sx}{z}\]
for all \(x \in \Ker f\) and all \(z \in Z\).
\end{proof}

\subsection{The canonical partial ordering and inner products}
\label{sect:r-rings}
Each of the examples of pre\nobreakdash-Hilbert \(\inv\)\nobreakdash-categories described in \cref{sect:fin-inner-prod-modules,sect:hilbert-W-modules} is a \(\inv\)\nobreakdash-subcategory of \(\InnerProd{R}\) for some ordered \(\inv\)\nobreakdash-ring~\(R\). The following proposition, which describes canonical \(\inv\)\nobreakdash-functors from a pre\nobreakdash-Hilbert \(\inv\)\nobreakdash-category to \(\InnerProd{R}\) for various \(R\), explains why this is no coincidence. A \textit{\(\inv\)\nobreakdash-functor} is a functor~\(F\) between \(\inv\)\nobreakdash-categories such that \(Ff^\inv = (Ff)^\inv\) for all morphisms \(f\).

\begin{proposition}
\label{prop:canonical-structure}
For all objects \(A\) of a pre\nobreakdash-Hilbert \(\inv\)\nobreakdash-category \(\C\),
\begin{itemize}
    \item the hom-set \(\C(A, A)\) is canonically an ordered \(\inv\)\nobreakdash-ring,
    and
    \item 
    the hom-functor \(\C(A, \blank) \colon \C \to \Set\) factors canonically through the forgetful functor \(\InnerProd{\C(A, A)} \to \Set\) via a \(\inv\)\nobreakdash-functor \(\C \to \InnerProd{\C(A, A)}\).
\end{itemize}
In particular, for all objects \(A\) and \(X\) of \(\C\),
\begin{itemize}
    \item addition in the ring \(\C(A, A)\) and in the module \(\C(A, X)\) both come from the unique enrichment of \(\C\) in abelian groups,
    \item multiplication in \(\C(A, A)\) and scalar multiplication in \(\C(A, X)\) are both just composition in \(\C\),
    \item the involution on \(\C(A, A)\) is the restriction of the involution on \(\C\),
    \item the positive cone on \(\C(A, A)\) is defined by
    \[\Pos{\C(A, A)} = \setb[\big]{y^\inv y}{y \in \C(A, Y) \text{ for some object } Y \text{ of } \C},\]
    \item the inner product on \(\C(A, X)\) is defined by \(\innerProd{x}{y} = x^\inv y\).
\end{itemize}
\end{proposition}

\begin{proof}[Proof of \cref{prop:canonical-structure}]
By \cref{cor:additive}, the category \(\C\) is additive, and thus uniquely enriched in abelian groups. In particular, the hom-sets \(\C(A, A)\) and \(\C(A, X)\) are canonically abelian groups, and composition distributes over the group operations. Hence \(\C(A, A)\) is a ring with addition given by the group operation and multiplication given by composition, and \(\C(A, X)\) is a right \(\C(A, A)\)\nobreakdash-module with addition given by the group operation and scalar multiplication given by composition.

The ring \(\C(A, A)\) is \(\inv\)\nobreakdash-ring when equipped with the restriction of the involution on \(\C\). Indeed, \((r + s)^\inv = r^\inv + s^\inv\) for all \(r, s \in \C(A, A)\) by \cref{prop:star-addition}, and the remaining \(\inv\)\nobreakdash-ring axioms come from the corresponding \(\inv\)\nobreakdash-category axioms. The \(\inv\)\nobreakdash-ring \(\C(A, A)\) is anisotropic by~\cite[Lemma~2.5]{vicary:2011:complex-numbers}. The set \(\Pos{\C(A, A)}\) is a positive cone. Indeed, for all \(x \in \C(A, X)\), \(y \in \C(A, Y)\) and \(r \in \C(A, A)\),
\begin{align*}
    (x^\inv x)^\inv &= x^\inv (x^{\inv})^\inv = x^\inv x,&
    1 &= 1^\inv 1,\\
    x^\inv x + y^\inv y &= \pairBig{x}{y}^\inv\pairBig{x}{y},&
    r^\inv x^\inv x r &= (xr)^\inv xr.
\end{align*}
Also, if \(x^\inv x = -y^\inv y\), then \(0 = \pair{x}{y}^\inv \pair{x}{y}\), so \(\pair{x}{y} = 0\) by~\cite[Lemma~2.5]{vicary:2011:complex-numbers}, and thus
\[x^\inv x = x^\inv p_1\pairBig{x}{y} = x^\inv p_1 0 = 0.\]

We now verify the inner product axioms. Linearity follows from associativity of composition and distributivity of composition over addition. Hermitianness follows from the \(\inv\)\nobreakdash-category axioms. Positivity holds by definition. Anisotropy follows again from~\cite[Lemma~2.5]{vicary:2011:complex-numbers}.

Let \(f \colon X \to Y\) be a morphism in \(\C\). For each \(x \in \C(A, X)\) and each \(y \in \C(A, Y)\),
\[\innerProd[\big]{y}{\C(A, f)x} = y^\inv \! fx = (f^\inv y)^\inv x = \innerProd[\big]{\C(A, f^\inv)y}{x}.\]
Hence the function \(\C(A, f)\) is adjointable with respect to the canonical inner product \(\C(A, A)\)\nobreakdash-module structures on \(\C(A, X)\) and \(\C(A, Y)\), and \(\C(A, f)^\inv = \C(A, f^\inv)\).
\end{proof}

The proof of \cref{prop:canonical-structure} appealed twice to the fact~\cite[Lemma~2.5]{vicary:2011:complex-numbers} that in \(\inv\)\nobreakdash-categories with a zero object and isometric equalisers, if \(f^\inv f = 0\) then \(f = 0\). For pre\nobreakdash-Hilbert \(\inv\)\nobreakdash-categories, this fact is now completely captured by the anisotropy of the canonical inner products on their hom-sets. Anisotropy is a useful tool for proving equality of morphisms that is not available in (ordinary) category theory.

The specialisation of the partial ordering defined in \cref{prop:canonical-structure} to \(\Hilb\) is the usual partial ordering of the Hermitian operators on a Hilbert space. Actually, this is true, more generally, of its specialisation to \(\Hilb_A\)  for each complex W*-algebra \(A\).

\begin{proposition}
\label{r:w-algebra-order}
Let \(A\) be a complex W*-algebra. Let \(f \colon X \to X\) be a morphism of \(\Hilb_A\). Then \(f \geq 0\) if and only if \(\innerProd{x}{fx} \geq 0\) for each \(x \in X\).
\end{proposition}

\begin{proof}
If \(f \geq 0\), then, by definition, there is a morphism \(g \colon X \to Y\) in \(\Hilb_A\) such that \(f = g^\inv g\), and so \(\innerProd{x}{fx} = \innerProd{x}{g^\inv gx} = \innerProd{gx}{gx} \geq 0\) for each \(x \in X\). Conversely, if \(\innerProd{x}{fx} \geq 0\) for each \(x \in X\), then~\cite[143~IV and 144~I]{westerbaan:dagger-and-dilation} there exists \(h \in \Hilb_A(X, X)\) such that \(f = h^\inv h\); in particular, \(f \geq 0\).
\end{proof}

\section{Positivity and inverses}
\label{sect:positivity-inverses}

If \(a\) and \(b\) are Hermitian elements of an ordered \(\inv\)\nobreakdash-ring \(R\), write \(a \slt b\) or \(b \sgt a\) if \(a \leq b\) and \(b - a\) is invertible. The element \(a\) is \textit{strictly positive} if \(a \sgt 0\). This section is about useful properties of strictly positive endomorphisms in pre\nobreakdash-Hilbert \(\inv\)\nobreakdash-categories.

First, inverses of strictly positive elements are strictly positive.

\begin{proposition}
\label{prop:inverse-strict-pos}
In an ordered \(\inv\)\nobreakdash-ring, if \(a \sgt 0\) then \(a^{-1} \sgt 0\).
\end{proposition}

\begin{proof}
The element \(a^{-1}\) has inverse \(a\). It is also Hermitian. Indeed, 
\[(a^{-1})^\inv a = (a^{-1})^\inv a^\inv = (aa^{-1})^\inv = 1^\inv = 1,\]
and similarly \(a(a^{-1})^\inv = 1\), so \(a^{-1} =(a^{-1})^\inv\). Finally,
\[a^{-1} = 1a^{-1} = (a^{-1})^\inv a a^{-1} \geq (a^{-1})^\inv 0 a^{-1} = 0. \qedhere\]
\end{proof}

Next, inversion of strictly positive endomorphisms is anti-monotone. The following proposition is a generalisation of this fact.

\begin{proposition}
\label{prop:inv_antitone}
For all morphisms 
\[
    \begin{tikzcd}[cramped]
        A
            \arrow[r, "f"]
            \arrow[loop, "a", looseness=10, in=155, out=-155, overlay]
            \&
        B   \arrow[loop, "b" swap, looseness=10, in=25, out=-25, overlay]
    \end{tikzcd}
\]
in a pre\nobreakdash-Hilbert \(\inv\)\nobreakdash-category, if \(a \sgt 0\) and \(b \sgt 0\) then
\[
    f^\inv b^{-1} f \leq a    
    \;\iff\;
    fa^{-1}\!f^\inv \leq b
    \qquad\text{and}\qquad
    f^\inv b^{-1} f \slt a
    \;\iff\;
    fa^{-1}\!f^\inv \slt b.
\]
\end{proposition}

\begin{proof}
The \textit{if} and \textit{only if} directions of the equivalences are equivalent to each other by swapping \(f\) with \(f^\inv\) and \(a\) with \(b\). If \(a - f^{\inv}b^{-1}\!f\) is positive, then
\[
    b - fa^{-1}\!f^\inv
    = fa^{-1}(a - f^\inv b^{-1}\!f)a^{-1}\!f^\inv
    + (1 - f a^{-1}\!f^{\inv}b^{-1})b(1 - b^{-1}\!f a^{-1}\!f^{\inv})
\]
is also positive. If \(a - f^{\inv}b^{-1}\!f\) is invertible, then
\[(b - fa^{-1}\!f^\inv)^{-1} = b^{-1} + b^{-1}\!f(a - f^{\inv}b^{-1}\!f)^{-1}\!f^\inv b^{-1}. \qedhere\]
\end{proof}

\cref{prop:inv_antitone} is closely related to \textit{Schur complements} and the \textit{Woodbury matrix identity}~\cite{zhang:schur-complement}. It will be used in \cref{sect:contractions} to show that (strict) contractions are preserved by the involution. A similar result holds for all elements \(a\), \(b\) and \(f\) of an ordered \(\inv\)\nobreakdash-ring.

Next, strictly positive endomorphisms are represented by isomorphisms.

\begin{proposition}
\label{p:strict-pos}
In a pre\nobreakdash-Hilbert \(\inv\)\nobreakdash-category, a morphism \(a \colon A \to A\) is strictly positive if and only if \(a = x^\inv x\) for some object \(X\) and some isomorphism \(x \colon A \to X\).
\end{proposition}

\begin{proof}
The \textit{if} direction is trivial. For the \textit{only if} direction, suppose that \(a \sgt 0\). As \(a \geq 0\), there is an object \(Y\) and a morphism \(y \colon A \to Y\) such that \(a = y^\inv y\). As \(y^\inv y = a\) is invertible, the morphism \(y\) is a closed monomorphism and \((y^\inv y)^{-1} y^\inv\) is its canonical retraction. In particular, \(y\) is a kernel of the idempotent \(1 - y(y^\inv y)^{-1} y^\inv\). By axiom \cref{axiom:R-kernel}, the morphism \(1 - y(y^\inv y)^{-1} y^\inv\) also has an isometric kernel \(m \colon X \to Y\), as depicted in \cref{f:strict-pos}.
\begin{diagram}
    \centering
    \begin{tikzcd}[column sep=large]
    A
        \arrow[dr, "y"]
        \arrow[d, "x" swap]
        \&
        \&[3em]
    \\
    X
        \arrow[r, "m" swap]
        \&
    Y
        \arrow[r, "1 - y(y^\inv y)^{-1} y^\inv"]
        \&
    Y
    \end{tikzcd}
    \caption{}
    \label{f:strict-pos}
\end{diagram}
Universality of the kernels \(y\) and \(m\) yields an isomorphism \(x \colon A \to X\) such that \(mx = y\). Observe that
\[a = y^\inv y = (mx)^\inv mx = x^\inv m^\inv m x = x^\inv 1 x = x^\inv x.\qedhere\]
\end{proof}

We now focus on a property of ordered \(\inv\)\nobreakdash-rings called \textit{symmetry}, which comes from unbounded operator algebra (see, e.g., \cite{schotz:equivalence-order-algebraic,schmudgen:2012:unbounded-operators,dixon:generalized-b-algebras,inoue:class-unbounded-operator}), and a stronger variant of symmetry that we call \textit{inverse closure}.

\begin{definition}
An ordered \(\inv\)\nobreakdash-ring is 
\begin{itemize}
    \item \textit{symmetric} if \(a \geq 1\) implies \(a \sgt 0\), and
    \item \textit{inverse closed} if \(a \geq b\) and \(b \sgt 0\) implies \(a \sgt 0\).
\end{itemize}
\end{definition}

\noindent Clearly if an ordered \(\inv\)\nobreakdash-ring is inverse closed then it is symmetric.

Surprisingly, if \(\C\) is a pre\nobreakdash-Hilbert \(\inv\)\nobreakdash-category, then each of the ordered \(\inv\)\nobreakdash-rings \(\C(A, A)\) is inverse closed. This fact is used in \cref{rem:weak-bounded-transform} to show that every morphism of a pre\nobreakdash-Hilbert \(\inv\)\nobreakdash-category factors into a contraction and the inverse of a contraction.

\begin{proposition}
\label{prop:symmetric}
For each object \(A\) of a pre\nobreakdash-Hilbert \(\inv\)\nobreakdash-category \(\C\), the ordered \(\inv\)\nobreakdash-ring \(\C(A, A)\) is inverse closed.
\end{proposition}

The proof of a similar result~\cite[Proposition~1.13]{handelman:1981:rings-with-involution-partially} about Rickart \(\inv\)\nobreakdash-rings can be adapted, replacing right annihilators with kernels. Instead, we give a simpler proof. 

\begin{proof}
Let \(a, b \in \C(A, A)\) and suppose that \(a \geq b\) and \(b \sgt 0\). Since \(a \geq b\), there is an object \(X\) and a morphism \(x \colon A \to X\) such that \(a - b = x^\inv x\). Also, by \cref{p:strict-pos}, there is an object \(Y\) and an isomorphism \(y \colon A \to Y\) such that \(b = y^\inv y\).

Refer to \cref{f:inv-closed}.
\begin{diagram}
    \centering
    \begin{tikzcd}[column sep=huge]
    B
        \arrow[d, "y^{-1}v" swap]
        \arrow[dr, "{\pair{u}{v}}"]
        \&
        \&
    \\
    A
        \arrow[r, "{\pair{x}{y}}" swap]
        \&
    X \oplus Y
        \arrow[r, "{\copair{-1}{xy^{-1}}}" swap]
        \&
    X
    \end{tikzcd}
    \caption{}
    \label{f:inv-closed}
\end{diagram}
Let us show that \(\pair{x}{y}\) is a kernel of \(\copair{-1}{xy^{-1}}\). First,
\[
    \copairBig{-1}{xy^{-1}}\pairBig{x}{y} = -x + xy^{-1}y
    = -x + x = 0.
\]
Next, let \(\pair{u}{v} \colon B \to X \oplus Y\) and suppose that \(\copair{-1}{xy^{-1}}\pair{u}{v} = 0\). For uniqueness, if \(\pair{u}{v} = \pair{x}{y}f\) for some morphism \(f \colon B \to A\), then \(v = yf\), so \(f = y^{-1}v\). For existence,
\[
    0
    = \copairBig{-1}{xy^{-1}}\pairBig{u}{v}
    = -u + xy^{-1}v,
\]
so \(u = xy^{-1}v\), and thus
\[
    \pairBig{u}{v}
    = \pairBig{xy^{-1}v}{v}
    = \pairBig{x}{y}y^{-1}v.
\]

We have shown that \(\pair{x}{y} \colon A \to X \oplus Y\) is a normal monomorphism. By \cref{lem:split-is-closed}, it is also a closed monomorphism. Hence
\[
    a = (a - b) + b = x^\inv x + y^\inv y = \pairBig{x}{y}^\inv \pairBig{x}{y}
\]
is positive and invertible.
\end{proof}

Every non-zero symmetric ordered \(\inv\)\nobreakdash-ring contains a copy of the ordered \(\inv\)\nobreakdash-ring \(\Rats\) of rational numbers. An \textit{ordered \(\inv\)\nobreakdash-algebra} over a commutative ordered \(\inv\)\nobreakdash-ring \(R\) is an ordered \(\inv\)\nobreakdash-ring~\(A\) equipped with a \textit{ordered \(\inv\)\nobreakdash-ring homomorphism} from \(R\) to the centre \(\centre_A\) of~\(A\), that is, a ring homomorphism \(\varphi \colon R \to \centre_A\) such that \(\varphi r^\inv = (\varphi r)^\inv\) for all \(r \in R\), and \(\varphi \Pos{R} \subseteq \Pos{A}\).

\begin{proposition}
\label{prop:rats-algebra}
If an ordered \(\inv\)\nobreakdash-ring is symmetric, then it is uniquely an ordered \(\inv\)\nobreakdash-algebra over \(\Rats\).
\end{proposition}

\begin{proof}
Let \(A\) be a symmetric ordered \(\inv\)\nobreakdash-ring.
The ring \(\Ints\) of integers is an initial object in the category of rings and ring homomorphisms. Let \(\varphi\) be the unique ring homomorphism \(\Ints \to \centre_A\). Then \(\varphi\) is an ordered \(\inv\)\nobreakdash-ring homomorphism. Indeed, the involution on \(A\) restricts to a ring homomorphism \(\blank^\inv \colon \centre_A \to \centre_A\) because \(\centre_A\) is commutative, so uniqueness of \(\varphi\) ensures that \((\varphi n)^\inv = \varphi n\) for all \(n \in \Ints\). Also, \(\varphi 0 = 0 \in \Pos{A}\), and \(\varphi n \in \Pos{A}\) implies \(\varphi(1 + n) = 1 + \varphi n \in \Pos{A}\) for each \(n \in \Pos{\Ints}\), so \(\varphi \Pos{\Ints} \subseteq \Pos{A}\) by induction.

Recall that \(\Rats\) is the field of fractions of \(\Ints\). We will show that there is a unique extension \(\psi \colon \Rats \to \centre_A\) of \(\varphi\) along the embedding \(\iota \colon \Ints \hookrightarrow \Rats\). For uniqueness, let \(\psi \colon \Rats \to \centre_A\) be a ring homomorphism that extends \(\varphi\). Then
\begin{equation}
    \label{e:rat-embed-unique}
    \psi \paren[\Big]{\frac{m}{n}}
    = \psi\paren[\big]{(\iota m)(\iota n)^{-1}}
    = (\psi \iota m)(\psi \iota n)^{-1}
    = (\varphi m)(\varphi n)^{-1}.
\end{equation}
for all \(m, n \in \Ints\) with \(n \neq 0\). Hence \(\psi\) is determined by \(\varphi\).

For existence, first observe, for each \(n \in \Ints\), that \(n \neq 0\) implies \(\varphi n\) is invertible. Indeed, if \(n \geq 1\), then
\[\varphi n = \varphi(n - 1) + 1 \geq 1,\]
so \(\varphi n\) is invertible by symmetry. Also, if \(n \leq -1\), then 
\[-n = (1 - n) - 1 \geq (1 - n) + n = 1,\]
 so \((\varphi n)^{-1} = -\paren[\big]{\varphi (-n)}^{-1}\).

Next observe, for all \(m,m',n,n' \in \Ints\) with \(n,n' \neq 0\), that \(mn' = m'n\) implies
\begin{multline*}
    (\varphi m)(\varphi n)^{-1}
    = (\varphi m)(\varphi n')(\varphi n')^{-1}(\varphi n)^{-1}
    = (\varphi mn')(\varphi n')^{-1}(\varphi n)^{-1}
    \\= (\varphi m'n)(\varphi n')^{-1}(\varphi n)^{-1}
    = (\varphi m')(\varphi n)(\varphi n')^{-1}(\varphi n)^{-1}
    = (\varphi m')(\varphi n')^{-1},
\end{multline*}
because \(\varphi\) is a ring homomorphism and \(\centre_A\) is commutative.

Equation \cref{e:rat-embed-unique} thus defines a function \(\psi \colon \Rats \to \centre_A\), and it is actually a ring homomorphism because \(\varphi\) is a ring homomorphism and \(\centre_A\) is commutative. Next, \((\psi q)^\inv = \psi q\) for all \(q \in \Rats\) by uniqueness of \(\psi\), again using the fact that \(\blank^\inv \colon \centre_A \to \centre_A\) is a ring homomorphism. Finally, \(\psi \Pos{\Rats} \subseteq \Pos{A}\). Indeed, let \(q \in \Pos{\Rats}\). Then
\(q = \frac{m}{n}\) for some \(m, n \in \Pos{\Ints}\) with \(n \neq 0\). Now \(\varphi n \sgt 0\) because \(\varphi \Pos{\Ints} \subseteq \Pos{A}\) and \(\varphi n\) is invertible. Hence \((\varphi n)^{-1} \sgt 0\) by \cref{prop:inverse-strict-pos}. Altogether, 
\[\psi q = \psi \paren[\Big]{\frac{m}{n}} = \psi\paren[\big]{\underbrace{n^{-1} + \dots + n^{-1}}_{m \text{ times}}} = \underbrace{(\varphi n)^{-1} + \dots + (\varphi n)^{-1}}_{m \text{ times}} \in \Pos{A}. \qedhere\]
\end{proof}
\begin{proposition}
\label{prop:rat-linear}
Every pre\nobreakdash-Hilbert \(\inv\)\nobreakdash-category is uniquely \(\Rats\)\nobreakdash-linear, that is, there is a unique way to equip its hom-sets with the structure of a vector space over \(\Rats\) such that composition is bilinear. Additionally, \((f \cdot q)^\inv = f^\inv \cdot q\) for all morphisms \(f\) and all \(q \in \Rats\).
\end{proposition}

\begin{proof}
Let \(\C\) be a pre\nobreakdash-Hilbert \(\inv\)\nobreakdash-category. Since it is additive (\cref{cor:additive}), it is uniquely enriched in abelian groups. Also, by \cref{prop:canonical-structure,prop:symmetric,prop:rats-algebra}, for each object \(X\), there is a unique ordered \(\inv\)\nobreakdash-ring homomorphism \(\psi_X \colon \Rats \to \centre_{\C(X, X)}\).

For each integer \(n \geq 1\), the endomorphisms \(\psi_X n\) are natural in \(X\). Indeed, for each morphism \(f \colon X \to Y\),
\[
    f(\psi_X n)
    = f \paren[\big]{\underbrace{1_X + \dots + 1_X}_{n \text{ times}}}
    = \underbrace{f + \dots + f}_{n \text{ times}}
    = \paren[\big]{\underbrace{1_Y + \dots + 1_Y}_{n \text{ times}}}f
    = (\psi_Y n)f.
\]
Since the inverses of the components of a natural isomorphism form a natural transformation, the endomorphisms \(\psi_X(n^{-1}) = (\psi_X n)^{-1}\) are also natural in \(X\). Using naturality of the identity and zero endomorphisms, and the fact that each \(\psi_X\) is a ring homomorphism, it is now easy to show that the endomorphisms \(\psi_X q\) are natural in \(X\) for each \(q \in \Rats\).

Let \(f \cdot q = f(\psi_X q)\) for each morphism \(f \colon X \to Y\) and each \(q \in \Rats\). Using the fact that each \(\psi_X\) is a ring homomorphism, and also distributivity of composition over addition, it is routine to verify that \(\cdot\) makes each of the abelian groups \(\C(X,Y)\) into a vector space over \(\Rats\). Using naturality of \(\psi_X q\) in \(X\) for each \(q \in \Rats\), it is straightforward to show that composition is bilinear, and that \((f \cdot q)^\inv = f^\inv \cdot q\).

It remains to prove uniqueness. Suppose that there is another \(\Rats\)\nobreakdash-linear structure on \(\C\), and denote its scalar multiplication by \(\odot\). Since \(\C\) has a unique enrichment in abelian groups, the enrichment in abelian groups underlying this new \(\Rats\)\nobreakdash-linear structure must be the canonical one. It is routine to verify, for each object \(X\), that the map \(q \mapsto 1_X \odot q\) is a ring homomorphism \(\Rats \to \centre_{\C(X, X)}\). By the uniqueness argument in the proof of \cref{prop:rats-algebra}, it follows that \(1_X \odot q = \psi_X q\) for all \(q \in \Rats\). Hence, for each morphism \(f \colon X \to Y\) and each \(q \in \Rats\),
\[
    f \odot q
    = (f1_X) \odot q
    = f(1_X \odot q)
    = f(\psi_X q)
    = f \cdot q. \qedhere
\]
\end{proof}
\section{Douglas' lemma}
\label{s:douglas}
Douglas' lemma~\cite[Theorem~1]{douglas:majorization-factorization-range} is a well-known and important result relating majorisation, factorisation, and range inclusion of operators on Hilbert spaces. Surprisingly, almost all of Douglas' lemma holds more generally in every pre\nobreakdash-Hilbert \(\inv\)\nobreakdash-category, as we shall see in \cref{p:doug-unbounded,p:doug-bounded} below.

A \textit{lift} of a morphism \(y \colon A \to Y\) through a morphism \(f \colon X \to Y\) is a morphism \(x \colon A \to X\) such that \(fx = y\).

\begin{theorem}
\label{p:doug-unbounded}
Morphisms \(f \colon X \to Y\) and \(y \colon A \to Y\) in a pre\nobreakdash-Hilbert \(\inv\)\nobreakdash-category satisfy \(\Ran y \leq \Ran f\) if and only if there is a lift of \(y\) through \(f\). Also, if such a lift exists, then there is a unique such lift \(x \colon A \to X\) that also satisfies \(\Ran x \leq (\Ker f)^\perp\). Furthermore, \(\Ker x = \Ker y\).
\end{theorem}

\begin{proof}
If there is a lift of \(y\) through \(f\), then \(\Ran y \leq \Ran f\) because \(\Ran y\) is the smallest subobject of \(Y\) through which \(y\) factors.

Conversely, suppose that \(\Ran y \leq \Ran f\). By \cref{s:range}, there is a coisometry \(e \colon X \to R\) and a monomorphism \(m \colon R \to Y\) such that \(f = me\). Since \(\Ran y \leq \Ran f\) and \(m\) represents \(\Ran f\), there is a morphism \(t \colon A \to R\) such that \(y = mt\). Let \(x = e^\inv t\). Then \(fx = mee^\inv t = mt = y\). Also, \(\Ran x \leq (\Ker f)^\perp\) because \(x\) factors through \(e^\inv\), which represents \((\Ker f)^\perp\).

To show that \(x\) is unique, consider a lift \(x' \colon A \to X\) of \(y\) through \(f\) such that \(\Ran x' \leq (\Ker f)^\perp\). Since \(e^\inv\) represents \((\Ker f)^\perp\), there is a morphism \(t' \colon A \to R\) such that \(x' = e^\inv t'\). Now \(m\) is monic and \(mt = y = fx' = mee^\inv t' = mt'\), so \(t = t'\), and thus \(x = e^\inv t = e^\inv t' = x'\).

Finally, \(\Ker x = \Ker e^\inv t = \Ker t = \Ker mt = \Ker y\) as \(e^\inv\) and \(m\) are monic.
\end{proof}

\begin{remark}
\label{r:semiclosed}
Consider the specialisation of \cref{p:doug-unbounded} to \(\Hilb\). It is easy to show that if \(\Ran y \leq \Ran f\) as subobjects of \(Y\), then \(\Ran y \subseteq \Ran f\) as subsets of \(Y\). To prove the converse, suppose that \(\Ran y \subseteq \Ran f\). Since \(\Ran y \leq \Ran f\) exactly when the inclusion function \(\Ran y \hookrightarrow \Ran f\) is continuous with respect to the norms on \(\Ran y\) and \(\Ran f\) respectively transferred from \(A\) and \(X\) as described in \cref{s:range}, it suffices, by the closed graph theorem, to show that the graph of this inclusion function is closed with respect to these norms. Consider a sequence \((v_n)_{n = 1}^\infty\) in \(\Ran y\) that converges to \(v \in \Ran y\) and \(u \in \Ran f\) with respect to the norms transferred from \(A\) and \(X\). Since the inclusion functions \(\Ran y \hookrightarrow Y\) and \(\Ran f \hookrightarrow Y\) are continuous with respect to these norms, the sequence also converges in \(Y\) to both \(u\) and \(v\), so \(u = v\), as required. By a similar argument, the partially ordered set of subobjects in \(\Hilb\) of a Hilbert space \(X\) is canonically isomorphic to the partially ordered set of semiclosed subspaces of \(X\); a subspace of \(X\) is \textit{semiclosed} (or an \textit{operator range}) if it can be equipped with a complete inner product that makes its inclusion into \(X\) continuous~\cite{fillmore:semiclosed}.
\end{remark}

An \textit{extension} of a morphism \(y \colon A \to Y\) along a morphism \(x \colon A \to X\) is a morphism \(f \colon X \to Y\) such that \(y = fx\).
The dual of \cref{p:doug-unbounded} says that if there is an extension of a morphism \(y \colon A \to Y\) along a morphism \(x \colon A \to X\), then there is a unique such extension \(f \colon X \to Y\) that also satisfies \(\Ker f \geq (\Ran x)^\perp\). This extension \(f\) will be called the \textit{reduced extension} of \(y\) along \(x\).

\begin{theorem}
\label{p:fact-eq}
Morphisms \(x \colon A \to X\) and \(y \colon A \to Y\) in a pre\nobreakdash-Hilbert \(\inv\)\nobreakdash-category satisfy \(x^\inv x = y^\inv y\) if and only if there is a partially isometric extension \(f\) of \(y\) along \(x\) such that \(\Ker f \leq (\Ran x)^\perp\). Moreover, if such an extension exists, then the reduced extension of \(y\) along \(x\) is also such an extension.
\end{theorem}

\begin{proof}
Suppose that there is a partial isometry \(f \colon X \to Y\) such that \(y = fx\) and \(\Ker f \leq (\Ran x)^\perp\). Then \(f(1 - f^\inv \! f) = 0\) as \(f\) is partially isometric. Hence \(1 - f^\inv \! f\) factors through \(\Ker f\), which in turn factors through~\((\Ran x)^\perp = \Ker x^\inv\). Thus \(x^\inv(1 - f^\inv \! f) = 0\), that is, \(x^\inv = x^\inv f^\inv \! f\). Hence \(y^\inv y = x^\inv f^\inv \! fx = x^\inv x\).

Conversely, suppose that \(x^\inv x = y^\inv y\). By the discussion at the end of \cref{sect:orthogonal-complements}, there is an isometry \(m \colon R \to X \oplus Y\) and an epimorphism \(e \colon A \to R\) such that \(me = \pair{x}{y}\). Let \(s = p_1m\) and \(t = p_2m\), so that \(m = \pair{s}{t}\).
Then \(x = se\) and \(y = te\). Now
\[
    e^\inv s^\inv se = x^\inv x = y^\inv y = e^\inv t^\inv te.
\]
Since \(e^\inv\) is monic and \(e\) is epic, it follows that \(s^\inv s = t^\inv t\). But \(1 = m^\inv m = s^\inv s + t^\inv t\), so \(s^\inv s = t^\inv t = 2^{-1}\). Let \(f = 2ts^\inv\). Then
\(ff^\inv \! f = 8ts^\inv st^\inv ts^\inv = 2ts^\inv = f\), so \(f\) is a partial isometry. Also, \(fx = 2ts^\inv se = te = y\). Finally,
\[
    \Ker f = \Ker 2ts^\inv = \Ker s^\inv = \Ker e^\inv s^\inv = \Ker x^\inv = (\Ran x)^\perp
\]
because \(e^\inv\) and \(2t\) are monic (\(2t\) is a section of \(t^\inv\)). In particular, \(f\) is the reduced extension of \(y\) along \(x\).
\end{proof}

\begin{corollary}
\label{p:fact-epic}
An epimorphism \(x \colon A \to X\) and a morphism \(y \colon A \to Y\) in a pre\nobreakdash-Hilbert \(\inv\)\nobreakdash-category satisfy \(x^\inv x = y^\inv y\) if and only if there is an isometric extension of \(y\) along \(x\). Moreover, if such an extension exists, then it is necessarily the reduced extension of \(y\) along \(x\).
\end{corollary}

\begin{proof}
This is a reduction of \cref{p:fact-eq} in the case where \(x\) is epic. Indeed, \((\Ran x)^\perp = \Ker x^\inv = 0\) because \(x^\inv\) is monic. Also, \(f\) is isometric if and only if it is partially isometric and monic, and it is monic if and only if \(\Ker f = 0\).
\end{proof}

\begin{theorem}
\label{p:doug-bounded}
Let \(q \in \Rats\). Morphisms \(x \colon A \to X\) and \(y \colon A \to Y\) in a pre\nobreakdash-Hilbert \(\inv\)\nobreakdash-category satisfy \(y^\inv y \leq q^2 x^\inv x\) if and only if there is an extension \(f\) of \(y\) along \(x\) such that \(f^\inv \! f \leq q^2\). Moreover, if such an extension exists, then the reduced extension of \(y\) along \(x\) is also such an extension.
\end{theorem}

\begin{proof}
If \(y = fx\) for some morphism \(f \colon X \to Y\) such that \(f^\inv \! f \leq q^2\), then
\[
    y^\inv y = x^\inv f^\inv \! fx \leq x^\inv q^2x = q^2x^\inv x.
\]

Conversely, suppose that \(y^\inv y \leq q^2x^\inv x\). Then there is a morphism \(z \colon A \to Z\) such that \(y^\inv y + z^\inv z = q^2x^\inv x\). Hence
\[
    \pairBig{y}{z}^\inv \pairBig{y}{z}
    = y^\inv y + z^\inv z
    = q^2x^\inv x
    = (xq)^\inv xq.
\]
By \cref{p:fact-eq}, there is a partial isometry \(w \colon X \to Y \oplus Z\) such that
\[
    wxq = \pairBig{y}{z}
    \qquad\text{and}\qquad
    \Ker w = (\Ran xq)^\perp,
\]
namely, the reduced extension of \(\pair{y}{z}\) along \(xq\). Let \(f = qp_1w\). By \cref{prop:con-is-subcategory}, the morphism \(p_1 w\) is a contraction. Hence \(f^\inv \! f = q (p_1w)^\inv p_1w q \leq q^2\). Also,
\[
    fx = qp_1wx = p_1wxq = p_1\pairBig{y}{z} = y.
\]
Finally, \(f\) is the reduced extension of \(y\) along \(x\). Indeed, if \(q = 0\), then
\[
    \Ker f = \Ker qp_1w = \Ker 0 = 1_X \geq (\Ran x)^\perp.
\]
If \(q \neq 0\), then \(q\) is monic, so
\[
    (\Ran xq)^\perp = \Ker qx^\inv = \Ker x^\inv = (\Ran x)^\perp,
\]
and thus
\[
    \Ker f = \Ker qp_1w \geq \Ker w = (\Ran xq)^\perp = (\Ran x)^\perp. \qedhere
\]
\end{proof}

\section{Contractions and codilators}
\label{sect:contractions}
This section develops the theory of contractions in pre\nobreakdash-Hilbert \(\inv\)\nobreakdash-categories. It is impressive how much of the theory of contractions on Hilbert spaces still holds in this general setting. The capstone of this section is a proof that each contraction in a pre\nobreakdash-Hilbert \(\inv\)\nobreakdash-category has a \textit{codilator}---a universal representation as a cospan of isometries. This is a variant of Sz.-Nagy's unitary dilation theorem~\cite{nagy:2010:harmonic-analysis}, which is the foundation for the modern theory of contractions on Hilbert spaces.

\subsection{Contractions}
\label{sect:category-con}
The notion of contraction in a pre\nobreakdash-Hilbert \(\inv\)\nobreakdash-category is defined in terms of the canonical partial ordering of the Hermitian morphisms, which was introduced in \cref{prop:canonical-structure}.

\begin{definition}
    \label{d:contraction}
A morphism \(f\) in a pre\nobreakdash-Hilbert \(\inv\)\nobreakdash-category is said to be \textit{contractive} if \(f^\inv \! f \leq 1\), and \textit{strictly contractive} if \(f^\inv \! f \slt 1\). A \textit{contraction} is a contractive morphism, and a \textit{strict contraction} is a strictly contractive morphism.
\end{definition}

The specialisation of \cref{d:contraction} to \(\Hilb\) is equivalent to the usual definition of contraction between Hilbert spaces. Actually, this is true more generally of its specialisation to \(\Hilb_A\) for each complex W*-algebra \(A\).

\begin{proposition}
    \label{p:hilb-module-norm-order}
Let \(A\) be a complex W*-algebra. A morphism \(f \colon X \to Y\) of \(\Hilb_A\) is contractive if and only if \(\norm{f} \leq 1\).
\end{proposition}

\begin{proof}
If \(f\) is contractive, then
\[\innerProd{fx}{fx} = \innerProd{x}{f^\inv \!fx} \leq \innerProd{x}{x}\]
for all \(x \in X\) by \cref{r:w-algebra-order}, and so~\cite[\S~2.9]{paschke:inner-product-modules-b-algebras}
\[\norm{f} = \inf \setb[\big]{c \in \PosReals}{\innerProd{fx}{fx} \leq c^2 \innerProd{x}{x}} \leq 1\]
for all \(x \in X\). Conversely, if \(\norm{f} \leq 1\), then~\cite[\S~2.6]{paschke:inner-product-modules-b-algebras} 
\[\innerProd{x}{f^\inv \!fx} = \innerProd{fx}{fx} \leq \norm{f}^2\innerProd{x}{x} \leq \innerProd{x}{x}\]
for all \(x \in X\), and so \(f^\inv \! f \leq 1\) by \cref{r:w-algebra-order}.
\end{proof}

\begin{proposition}
\label{prop:con-is-subcategory}
In a pre\nobreakdash-Hilbert \(\inv\)\nobreakdash-category,
\begin{enumerate}
    \item 
    \label{item:partial-isometry-contraction}
    every partial isometry is contractive,
    \item
    \label{item:composite-contraction}
    if \(f \colon X \to Y\) and \(g \colon Y \to Z\) are contractive, then \(gf\) is also contractive,
    \item
    \label{item:adjoint-contraction}
    if \(f\) is (strictly) contractive, then \(f^\inv\) is (strictly) contractive.
\end{enumerate}
\end{proposition}

\begin{proof}
If \(f\) is partially isometric, then
\[
    f^\inv \! f
    \leq f^\inv \! f + (1 - f^\inv \! f)^\inv(1 - f^\inv \! f)
    = 1,
\]
so property \cref{item:partial-isometry-contraction} holds. For property \cref{item:composite-contraction}, the morphism \(gf\) is contractive because
\[(gf)^\inv gf = f^\inv g^\inv g f \leq f^\inv \! f \leq 1\]
by the axioms of ordered \(\inv\)\nobreakdash-rings. Property \cref{item:adjoint-contraction} is the special case of \cref{prop:inv_antitone} where \(a = b = 1\).
\end{proof}

The contractions in a pre\nobreakdash-Hilbert \(\inv\)\nobreakdash-category \(\C\) thus form a wide \(\inv\)\nobreakdash-subcategory, which will be denoted \(\Con{\C}\). The following proposition lists several of its properties. Many of them play an important role in \citeauthor{heunenkornellvanderschaaf:con}'s characterisation of the category \(\Con{\Hilb}\) of Hilbert spaces and contractions~\cite{heunenkornellvanderschaaf:con}.

\begin{proposition}
\label{prop:contraction-properties}
Let \(\C\) be a pre\nobreakdash-Hilbert \(\inv\)\nobreakdash-category.
\begin{enumerate}
    \item 
    \label{prop:unitaries}
    Every split monomorphism in \(\Con{\C}\) is isometric.
    \item
    \label{prop:preserves-monicity}
    The inclusion functor \(\Con{\C} \hookrightarrow \C\) preserves jointly monic wide spans.
    \item 
    \label{prop:creates-isometric-equalisers}
    The inclusion functor \(\Con{\C} \hookrightarrow \C\) creates isometric equalisers. In particular, in \(\Con{\C}\), every parallel pair has an isometric equaliser.
    \item 
    \label{prop:isometric-equalisers}
    In \(\Con{\C}\), every equaliser is isometric.
\end{enumerate}
\end{proposition}

\begin{proof}
For property \cref{prop:unitaries}, if \(s \colon X \to Y\) is a split monomorphism in \(\Con{\C}\), then it has a retraction \(r \colon Y \to X\) in \(\Con{\C}\), and
\[s^\inv \! s \leq 1 = (r s)^\inv rs = s^\inv r^\inv r s \leq s^\inv \! s\]
so \(s^\inv \! s = 1\) by antisymmetry.

For property \cref{prop:preserves-monicity}, let \((g_\alpha \colon Y \to Z_\alpha)_{\alpha \in I}\) be a jointly monic wide span in \(\Con{\C}\), and let \(f_1, f_2 \colon X \to Y\) be morphisms in \(\C\) such that \(g_\alpha f_1 = g_\alpha f_2\) for all \(\alpha \in I\). As
\[
    f_1(1 + {f_1}^{\!\inv} \! f_1 + {f_2}^{\!\inv} \! f_2)^{-1}
    \qquad\text{and}\qquad
    f_2(1 + {f_1}^{\!\inv} \! f_1 + {f_2}^{\!\inv} \! f_2)^{-1}
\]
are contractions (see \cref{rem:weak-bounded-transform} below), and
\[g_\alpha f_1(1 + {f_1}^{\!\inv} \! f_1 + {f_2}^{\!\inv} \! f_2)^{-1} = g_\alpha f_2(1 + {f_1}^{\!\inv} \! f_1 + {f_2}^{\!\inv} \! f_2)^{-1}\]
for all \(\alpha \in I\), joint monicity of \((g_\alpha \colon Y \to Z_\alpha)_{\alpha \in I}\) in \(\Con{\C}\) implies that
\[f_1(1 + {f_1}^{\!\inv} \! f_1 + {f_2}^{\!\inv} \! f_2)^{-1} = f_2(1 + {f_1}^{\!\inv} \! f_1 + {f_2}^{\!\inv} \! f_2)^{-1}.\]
Hence \(f_1 = f_2\).

For property \cref{prop:creates-isometric-equalisers}, let \(f_1, f_2 \colon X \to Y\) be morphisms in \(\Con{\C}\), and let \(s \colon A \to X\) be an isometric equaliser of them in \(\C\). Then \(s\) comes from \(\Con{\C}\) because it is isometric. Also, if \(t \colon B \to X\) is a morphism in \(\Con{\C}\) such that \(f_1t = f_2t\), then there is a unique morphism \(t' \colon B \to A\) in \(\C\) such that \(t = st'\). But \(t'\) is also from \(\Con{\C}\) because \(t' = s^\inv \! s t' = s^\inv t\). Hence \(s\) is also an equaliser of \(f_1\) and \(f_2\) in \(\Con{\C}\).

For property \cref{prop:isometric-equalisers}, let \(f_1, f_2 \colon X \to Y\) be morphisms in \(\Con{\C}\), and let \(s \colon A \to X\) be an equaliser of them in \(\Con{\C}\). By property \cref{prop:creates-isometric-equalisers}, the morphisms \(f_1\) and \(f_2\) also have an isometric equaliser \(s' \colon A' \to X\) in \(\Con{\C}\). Let \(u\) be the unique morphism in \(\Con{\C}\) such that \(s = s'u\). Then \(u\) is an isomorphism in \(\Con{\C}\), so, by property \cref{prop:unitaries}, it is actually unitary. As \(s = s'u\) is a composite of isometries, it is also an isometry.
\end{proof}

\begin{remark}
    \label{rem:weak-bounded-transform}
    The proof of property \cref{prop:preserves-monicity} exhibits a useful technique for expressing morphisms \(f\) of a pre\nobreakdash-Hilbert \(\inv\)\nobreakdash-category in terms of contractions. Suppose that \(a \geq f^\inv \! f\). Then \(1 + a\) is invertible by \cref{prop:symmetric}. Hence
    \(f = f(1 + a)^{-1}(1 + a)\), where \(f(1 + a)^{-1}\) is contractive because
    \[
        (1 + a)^{-1}\! f^\inv \! f (1 + a)^{-1} \leq (1 + a)^{-1} (1 + a) (1 + a)^{-1} = (1 + a)^{-1} \leq 1
    \]
    by \cref{prop:inv_antitone}, and \((1 + a)^{-1}\) is similarly also contractive. This idea is inspired by the \textit{bounded transform} from unbounded operator algebra~(see, e.g., \cite{kaufman:1978:closed-operator-quotient,kaufman:closed-operators-contractions,koliha:kaufmans-theorem}).
\end{remark}

\subsection{Codilators}
\label{s:codilator}
Recall the following definition from the introduction.

\begin{definition}
A \textit{codilation} of a morphism \(f \colon X \to Y\) in a \(\inv\)\nobreakdash-category is a cospan \((T, t_1, t_2)\) where \(t_1 \colon X \to T\) and \(t_2 \colon Y \to T\) are isometries and \({t_2}^{\!\inv} t_1 = f\). A \textit{codilator} of \(f\) is a codilation \((S, s_1, s_2)\) of \(f\) such that, for all codilations \((T, t_1, t_2)\) of \(f\), there is a unique isometry \(t \colon S \to T\) such that \(t_1 = ts_1\) and \(t_2 = ts_2\).
\end{definition}

The dual notions of \textit{dilation} and \textit{dilator} are defined dually in terms of spans of coisometries. If \((S, s_1, s_2)\) is a codilator of~\(f\), then \((S, {s_1}^{\!\inv}\!, {s_2}^{\!\inv})\) is a dilator of~\(f\). This means that the notions of dilator and codilator are equivalent to each other.

\begin{lemma}
\label{p:codil-pairing}
In a pre\nobreakdash-Hilbert \(\inv\)\nobreakdash-category, a cospan 
\[
    \begin{tikzcd}[cramped]
        X
            \arrow[r, "s_1"]
            \&
        S
            \&
        Y
            \arrow[l, "s_2" swap]
    \end{tikzcd}
\]
is a codilation of a morphism \(f \colon X \to Y\) if and only if
\[
    \copairBig{s_1}{s_2}^\inv \copairBig{s_1}{s_2} = \begin{bmatrix} 1 & f^\inv \\ f & 1 \end{bmatrix}.
\]
\end{lemma}

\begin{proof}
The equations obtained by identifying the corresponding entries of the matrices on the left-hand and right-hand sides of the equation above are
\[
    {s_1}^\inv s_1 = 1,
    \qquad
    {s_2}^\inv s_1 = f,
    \qquad\text{and}\qquad
    {s_2}^\inv s_2 = 1.
\]
These equations are precisely what it means for \((S, s_1, s_2)\) to be a codilation of \(f\).
\end{proof}

\begin{lemma}
\label{p:codil-joint-epic}
In a pre\nobreakdash-Hilbert \(\inv\)\nobreakdash-category, a codilation of a morphism \(f \colon X \to Y\) is a codilator of \(f\) if and only if it is jointly epic.
\end{lemma}

\begin{proof}
Let \((S, s_1, s_2)\) be a codilation of \(f\). Suppose first that it is jointly epic, and let \((T, t_1, t_2)\) be another codilation of~\(f\). By \cref{p:codil-pairing},
\[
    \copairBig{t_1}{t_2}^\inv \copairBig{t_1}{t_2}
    = \copairBig{s_1}{s_2}^\inv \copairBig{s_1}{s_2}.
\]
Since \(s_1\) and \(s_2\) are jointly epic, \(\copair{s_1}{s_2}\) is epic. Hence, by \cref{p:fact-epic}, there is a unique isometric extension \(t \colon S \to T\) of \(\copair{t_1}{t_2}\) along \(\copair{s_1}{s_2}\). Also, a morphism \(t \colon S \to T\) is such an extension if and only if \(ts_1 = t_1\) and \(ts_2 = t_2\). Hence \((S, s_1, s_2)\) is a codilator of \(f\).

Conversely, suppose that \((S, s_1, s_2)\) is a codilator of \(f\). By the discussion at the end of \cref{sect:orthogonal-complements}, there is an epimorphism \(e \colon X \oplus Y \to E\) and an isometry \(m \colon E \to S\) such that \(me = \copair{s_1}{s_2}\). Let \(e_1 = ei_1\) and \(e_2 = ei_2\), so that \(e = \copair{e_1}{e_2}\). The cospan \((E, e_1, e_2)\) is a codilation of \(f\) by \cref{p:codil-pairing} because
\[
    \copairBig{e_1}{e_2}^\inv \copairBig{e_1}{e_2}
    = e^\inv e
    = e^\inv m^\inv me
    = \copairBig{s_1}{s_2}^\inv \copairBig{s_1}{s_2}.
\]
By universality of the codilator \((S, s_1, s_2)\), there is a unique isometry \(u \colon S \to E\) such that \(us_1 = e_1\) and \(us_2 = e_2\). Then \(ume_1 = us_1 = e_1\) and similarly \(ume_2 = e_2\). Since \(e_1\) and \(e_2\) are jointly epic, it follows that \(um = 1\). Hence \(m = u^\inv u m = u^\inv\), so \(m\) is unitary. Therefore \(\copair{s_1}{s_2} = m\copair{e_1}{e_2}\) is epic, and thus \((S, s_1, s_2)\) is jointly epic.
\end{proof}
\begin{theorem}
\label{p:codilation-exist}
The following statements about a morphism \(f \colon X \to Y\) in a pre\nobreakdash-Hilbert \(\inv\)\nobreakdash-category are equivalent:
\begin{enumerate}
    \item \(f\) is contractive,
    \item \(f\) has a codilation,
    \item \(f\) has a codilator,
    \item \(\begin{bmatrix} 1 & f^\inv \\ f & 1 \end{bmatrix} \geq 0\).
\end{enumerate}  
\end{theorem}

\begin{proof}
Clearly (iii) implies (ii). Also (ii) implies (i) by \cref{prop:con-is-subcategory}. Furthermore (ii) and (iv) are equivalent by \cref{p:codil-pairing}. It remains to show that (i) implies (iii).

Suppose that \(f\) is contractive. Since \(1 - f^\inv \! f \geq 0\), there is a morphism \(g \colon X \to Z\) such that \(1 - f^\inv \! f = g^\inv g\). By the discussion at the end of \cref{sect:orthogonal-complements}, there is an epimorphism \(e \colon X \to D\) and an isometry \(m \colon D \to Z\) such that \(g = me\).

The cospan
\[
    \begin{tikzcd}[cramped]
        X
            \arrow[r, "\pair{f}{e}"]
            \&
        Y \oplus D
            \&
        Y
            \arrow[l, "i_1" swap]
    \end{tikzcd}
\]
is a codilation of \(f\). Indeed,
\[
    {i_1}^\inv \pairBig{f}{e} = p_1 \pairBig{f}{e} = f,
\]
the morphism \(i_1\) is isometric by definition, and \(\pair{f}{e}\) is isometric because
\[
    \pairBig{f}{e}^\inv \pairBig{f}{e} = f^\inv \! f + e^\inv e = f^\inv \! f + e^\inv m^\inv m e = f^\inv \! f + g^\inv g = 1.
\]

The cospan is also jointly epic. Indeed, suppose that
\[
    \copairBig{s_1}{s_2}\pairBig{f}{e} = \copairBig{t_1}{t_2}\pairBig{f}{e}
    \qquad\text{and}\qquad
    \copairBig{s_1}{s_2}i_1 = \copairBig{t_1}{t_2}i_1,
\]
that is, suppose that
\[
    s_1f + s_2e = t_1f + t_2e
    \qquad\text{and}\qquad
    s_1 = t_1.
\]
Then \(s_2e = t_2e\), so \(s_2 = t_2\) since \(e\) is epic. But this means that \(\copair{s_1\vphantom{t_1}}{s_2\vphantom{t_2}} = \copair{t_1}{t_2}\).

The cospan is thus a codilator of \(f\) by \cref{p:codil-joint-epic}.
\end{proof}

\begin{proposition}
The following statements about a morphism \(f \colon X \to Y\) in a pre\nobreakdash-Hilbert \(\inv\)\nobreakdash-category are equivalent:
\begin{enumerate}
    \item \(f\) is strictly contractive,
    \item \(f\) has a codilation \((S, s_1, s_2)\) that is a coproduct,
    \item \(\begin{bmatrix} 1 & f^\inv \\ f & 1 \end{bmatrix} \sgt 0\).
\end{enumerate}
Furthermore, every codilation \((S, s_1, s_2)\) of \(f\) that is a coproduct is a codilator of \(f\).
\end{proposition}

\begin{proof}
Suppose first that (i) holds. The morphism
\[
    \begin{bmatrix} 1 & f^\inv \\ f & 1 \end{bmatrix}
\]
is positive by \cref{p:codilation-exist}. It is easy to show that
\[
    \begin{bmatrix}
        (1 - f^\inv \! f)^{-1}
        &
        -(1 - f^\inv \! f)^{-1}\!f^\inv
        \\
        -f(1 - f^\inv \! f)^{-1}
        &
        1 + f(1 - f^\inv \! f)^{-1} f^\inv
    \end{bmatrix}
\]
is its inverse. Hence (iii) holds.

Conversely, suppose that (iii) holds. By \cref{p:codilation-exist}, \(f\) is contractive.
Letting
\[
    \begin{bmatrix}a & b\\ c & d \end{bmatrix}
    = \begin{bmatrix}1 & f^\inv \\ f & 1 \end{bmatrix}^{-1},
\]
it is easy to show that \(a\) is an inverse of \(1 - f^\inv \! f\). Hence (i) holds.

Suppose now that (iii) holds. By \cref{p:codilation-exist}, \(f\) has a codilator \((S, s_1, s_2)\). The morphism \(\copair{s_1}{s_2}\) is closed monic by \cref{p:codil-pairing}, and it is epic by \cref{p:codil-joint-epic}. Hence it is invertible. Thus (ii) holds.

Conversely, suppose that (ii) holds. Then \(\copair{s_1}{s_2}\) is invertible, so (iii) follows from \cref{p:codil-pairing}. Also, since \((S, s_1, s_2)\) is a jointly epic codilation of~\(f\), it is a codilator of \(f\) by \cref{p:codil-joint-epic}.
\end{proof}
The universal property of codilators resembles that of pushouts. The following proposition explains why this is no coincidence. For background on orthogonal complements in pre\nobreakdash-Hilbert \(\inv\)\nobreakdash-categories, see \cref{sect:orthogonal-complements}.

\begin{proposition}
\label{prop:codilator-pushout}
Consider isometries \(s \colon A \to X\) and \(t \colon A \to Y\) in a pre\nobreakdash-Hilbert \(\inv\)\nobreakdash-category. The square
\begin{equation}
    \label{eq:pushout-square}
    \begin{tikzcd}[sep={3.0em,between origins}, column sep={3.6em,between origins}]
        \&
    X
        \arrow[dr, "{\bsmallmat{s^{\perp \inv}\\s^\inv\\0}}"]
        \&
    \\
    A
        \arrow[ur, "s"]
        \arrow[dr, "t" swap]
        \&
        \&
    (X \ominus A) \oplus A \oplus \mathrlap{(Y \ominus A)}
    \\
        \&
    Y
        \arrow[ur, "{\bsmallmat{0\\t^\inv\\t^{\perp \inv}}}" swap]
    \end{tikzcd}
\end{equation}
is a pushout square. Also, the cospan
\begin{equation}
\label{eq:pushout-codilator}
    \begin{tikzcd}[cramped, sep=huge]
        X
            \arrow[r, "{\bsmallmat{s^{\perp \inv}\\s^\inv\\0}}"]
        \&
        (X \ominus A) \oplus A \oplus (Y \ominus A)
        \&
        Y
            \arrow[l, "{\bsmallmat{0\\t^\inv\\t^{\perp \inv}}}" swap]
    \end{tikzcd}
\end{equation}
is a codilator of the partial isometry \(ts^\inv\).
\end{proposition}

\begin{proof}
We begin by showing that the square \cref{eq:pushout-square} is a pushout. It commutes because
\[
    \begin{bmatrix}
        s^{\perp\inv}\\
        s^\inv\\
        0
    \end{bmatrix}
    s
    = \begin{bmatrix}
        s^{\perp\inv} \! s\\
        s^\inv \! s\\
        0s
    \end{bmatrix}
    =
    \begin{bmatrix}
        0\\
        1\\
        0
    \end{bmatrix}
    = \begin{bmatrix}
        0t\\
        t^\inv t\\
        t^{\perp\inv} t
    \end{bmatrix}
    = \begin{bmatrix}
        0\\
        t^\inv\\
        t^{\perp\inv}
    \end{bmatrix}
    t.
\]
Let \(f \colon X \to Z\) and \(g \colon Y \to Z\) be morphisms such that \(fs = gt\). For uniqueness, suppose that there is a morphism \(h \colon (X \ominus A) \oplus A \oplus (Y \ominus A) \to Z\) such that
\begin{equation}
\label{eq:comparison-morphism}
    h\begin{bmatrix}
        s^{\perp\inv}\\
        s^\inv\\
        0
    \end{bmatrix} = f
    \qquad\text{and}\qquad 
    h\begin{bmatrix}
        0\\
        t^\inv\\
        t^{\perp\inv}
    \end{bmatrix} = g.
\end{equation}
Then
\[
    hi_2
    =
    h\begin{bmatrix}
        0\\
        1\\
        0
    \end{bmatrix}
    =
    h\begin{bmatrix}
        s^{\perp\inv} \! s\\
        s^\inv \! s\\
        0s
    \end{bmatrix}
    = 
    h\begin{bmatrix}
        s^{\perp\inv}\\
        s^\inv\\
        0
    \end{bmatrix}s
    =
    fs = gt.
\]
Also
\[
    hi_1
    =
    h\begin{bmatrix}
        1\\
        0\\
        0
    \end{bmatrix}
    =
    h\begin{bmatrix}
        s^{\perp\inv} \! s^\perp\\
        s^\inv \! s^\perp\\
        0s^\perp
    \end{bmatrix}
    =
    h\begin{bmatrix}
        s^{\perp\inv}\\
        s^\inv\\
        0
    \end{bmatrix}s^\perp
    =
    fs^\perp,
\]
and \(hi_3 = gt^\perp\) similarly. Hence
\begin{equation}
\label{eq:pushout-uniqueness}
    h
    = \begin{bmatrix}fs^\perp & fs & gt^\perp\end{bmatrix}
    = \begin{bmatrix}fs^\perp & gt & gt^\perp\end{bmatrix},
\end{equation}
and this equation determines \(h\). For existence, if \(h\) is defined using equation \cref{eq:pushout-uniqueness}, then
\begin{multline*}
    h
    \begin{bmatrix}
        s^{\perp\inv}\\
        s^\inv\\
        0
    \end{bmatrix}
    =
    \begin{bmatrix}fs^\perp & fs & gt^\perp\end{bmatrix}
    \begin{bmatrix}
        s^{\perp\inv}\\
        s^\inv\\
        0
    \end{bmatrix}
    \\=
    fs^\perp s^{\perp\inv} + fss^\inv + gt^\perp0
    = f(s^\perp s^{\perp\inv} + ss^\inv) + 0
    = f1
    = f
\end{multline*}
and similarly \(h \bsmallmat{
    0\\
    t^\inv\\
    t^{\perp\inv}\\
} = g\). Hence the square \cref{eq:pushout-square} is a pushout.

We now show that the cospan \cref{eq:pushout-codilator} is a codilator of the morphism \(ts^\inv\). First,
\[
    \begin{bmatrix}
        s^{\perp\inv}\\
        s^\inv\\
        0
    \end{bmatrix}^\inv    
    \begin{bmatrix}
        s^{\perp\inv}\\
        s^\inv\\
        0
    \end{bmatrix}
    =
    \begin{bmatrix}
        s^{\perp} &
        s &
        0
    \end{bmatrix}
    \begin{bmatrix}
        s^{\perp\inv}\\
        s^\inv\\
        0
    \end{bmatrix}
    = s^{\perp}s^{\perp\inv} + ss^\inv + 0
    = 1
\]
and \(
\bsmallmat{
    0\\
    t^\inv\\
    t^{\perp\inv}\\
}^\inv \bsmallmat{
    0\\
    t^\inv\\
    t^{\perp\inv}\\
} = 1
\) similarly. Also
\[
    \begin{bmatrix}
        0\\
        t^\inv\\
        t^{\perp\inv}\\
    \end{bmatrix}^\inv
    \begin{bmatrix}
        s^{\perp\inv}\\
        s^\inv\\
        0
    \end{bmatrix}
    = 
    \begin{bmatrix}
        0 &
        t &
        t^\perp
    \end{bmatrix}
    \begin{bmatrix}
        s^{\perp\inv}\\
        s^\inv\\
        0
    \end{bmatrix}
    = 0s^{\perp\inv} + ts^\inv + t^\perp0
    = ts^\inv.
\]
Hence it is indeed a codilation of \(ts^\inv\). Let \((Z, f, g)\) be another codilation of \(ts^\inv\). Then
\begin{align*}
    (fs - gt)^\inv (fs - gt)
    &= s^\inv \!f^\inv \! f s - t^\inv \!g^\inv \! fs - s^\inv\! f^\inv\! gt + t^\inv \! g^\inv \! g t
    \\
    &= s^\inv \!f^\inv \! f s - t^\inv ts^\inv \!s - s^\inv \!st^\inv t + t^\inv \! g^\inv \! g t
    \\&= 1 - 1 - 1 + 1
    \\&= 0,
\end{align*}
because \(f\), \(g\), \(s\) and \(t\) are isometries and \(f^\inv \!g = st^\inv\). As the canonical inner product on \(\C(A,Z)\) is anisotropic, it follows that \(fs = gt\). But \cref{eq:pushout-square} is a pushout square. Hence the morphism \(h\) defined by equation \cref{eq:pushout-uniqueness} is the unique morphism \(h\) that satisfies equations \cref{eq:comparison-morphism}. It is isometric because
\begin{align*}
    h^\inv h
    &=
    \begin{bmatrix}s^{\perp\inv}\!f^\inv\\ s^\inv \!f^\inv \\ t^{\perp\inv}\!g^\inv\end{bmatrix}
    \begin{bmatrix}fs^\perp & fs & gt^\perp\end{bmatrix}
    =
    \begin{bmatrix}
        s^{\perp\inv}\!f^\inv \!fs^\perp
        &
        s^{\perp\inv}\!f^\inv \! fs
        &
        s^{\perp\inv}\!f^\inv \!gt^\perp
    \\
        s^\inv \!f^\inv \!fs^\perp
        &
        s^\inv \!f^\inv \! fs
        &
        t^\inv \!g^\inv \!gt^\perp
    \\
        t^{\perp\inv}\!g^\inv \! fs^\perp
        &
        t^{\perp\inv}\!g^\inv \! gt
        &
        t^{\perp\inv}\!g^\inv \! gt^\perp
    \end{bmatrix}
    \\&\qquad\qquad\qquad\qquad=
    \begin{bmatrix}
        s^{\perp\inv}1s^\perp
        &
        s^{\perp\inv}1s
        &
        s^{\perp\inv}\!st^\inv t^\perp
    \\
        s^\inv 1s^\perp
        &
        s^\inv 1s
        &
        t^\inv 1t^\perp
    \\
        t^{\perp\inv}ts^\inv \!s^\perp
        &
        t^{\perp\inv}1t
        &
        t^{\perp\inv}1t^\perp
    \end{bmatrix}
    =
    \begin{bmatrix}
        1
        &
        0
        &
        0
    \\
        0
        &
        1
        &
        0
    \\
        0
        &
        0
        &
        1
    \end{bmatrix}
    = 1. \qedhere
\end{align*}
\end{proof}

\noindent The following two degenerate cases of \cref{prop:codilator-pushout} are also notable:
\begin{itemize}
    \item for all objects \(X\) and \(Y\), the orthonormal coproduct \((X \oplus Y,\,i_1, i_2)\) is a codilator of the zero morphism \(0 \colon X \to Y\); and
    \item for each isometry \(t \colon A \to Y\), the cospan \((Y, t, 1)\) is a codilator of \(t\).
\end{itemize}
They are obtained by setting \(A = \zero\) and \(s = 1\), respectively.
\newpage
\printbibliography
\end{document}